\documentclass[11pt]{amsart}%
\usepackage{amsfonts}
\usepackage{amssymb,amsmath,amsthm}
\usepackage{amsmath}
\usepackage{amssymb}
\usepackage{graphicx}%
\providecommand{\U}[1]{\protect\rule{.1in}{.1in}}
\newtheorem{theorem}{Theorem}[section]
\newtheorem{lemma}[theorem]{Lemma}
\newtheorem{proposition}[theorem]{Proposition}
\newtheorem{corollary}[theorem]{Corollary}
\theoremstyle{definition}
\newtheorem{definition}[theorem]{Definition}
\newtheorem{example}[theorem]{Example}
\theoremstyle{remark}
\newtheorem{remark}[theorem]{Remark}

\makeatletter
\@namedef{subjclassname@2020}{\textup{2020} Mathematics Subject Classification}
\makeatother
\begin{document}
\title[Simple modules generated by indicator functions]{Simple modules over Leavitt path algebras\\generated by indicator functions on paths}
\author{Elizabeth Pacheco}
\address{School of Computer, Data and Mathematical Sciences, Western Sydney University,
Sydney, Australia}
\email{e.pacheco@westernsydney.edu.au}
\author{Kulumani M. Rangaswamy}
\address{Department of Mathematics, University of Colorado, Colorado Springs, Colorado
80918, USA}
\email{krangasw@uccs.edu}

\begin{abstract}
Let $E$ be an arbitrary graph, $K$ a field, and $D$ the set of boundary paths
of $E$ consisting of both infinite paths and finite paths ending at singular
vertices. The space $K^{D}$ of $K$-valued functions on $D$ is a right module
over the Leavitt path algebra $L_{K}(E)$. For a non-empty set $S\subseteq D$
that is strongly connected and closed under adding and deleting initial
segments, we study the submodule $N_{S}=\mathbf{1}_{S}\cdot L_{K}(E)$
generated by the indicator function of $S$. Let $S^{0}$ be the set of vertices
lying on the paths in $S$. We show that $N_{S}$ is simple if and only if $S$
contains no breaking vertex of the hereditary saturated set $H=E^{0}\setminus
S^{0}$. In particular, when $S$ contains only infinite paths or when $E$ is a
row-finite graph, $N_{S}$ is a simple module. Up to isomorphism, the simple
module $N_{S}$ depends only on $S^{0}$. When $S$ is an infinite set of paths,
the simple module $N_{S}$ is neither graded nor projective. The left module
$L_{K}(E)\mathbf{1}_{S}$, obtained by twisting the action through the
involution of $L_{K}(E)$, is simple exactly when $N_{S}$ is, and, when $S$ is
infinite, it is not isomorphic to any of the Chen simple modules. For the rose
graph $R_{n}$ with $n\geq2$ petals, we show that $L_{K}(R_{n})\mathbf{1}_{S}$
is isomorphic to one of the twisted Chen modules constructed by Nam,
Srivastava and Vien, and to none of the other simple modules in their list.
Moreover, when $S$ is not a single vertex, the simple module $N_{S}$ is
finitely presented if and only if only finitely many edges of $E$ have both
source and range in $S^{0}$; in particular, every simple module $N_{S}$ is
finitely presented when $E$ is a finite graph. Finally, we determine when the
analogous module generated by an indicator function over the Cohn path algebra
$C_{K}(E)$ is simple.

\end{abstract}
\maketitle

\section{Introduction}

Let $E$ be an arbitrary graph and $K$ a field. Simple modules over the Leavitt
path algebra $L_{K}(E)$ have been studied in a series of papers, both as a
first step towards the module theory of $L_{K}(E)$ and for their impact on the
algebraic properties of $L_{K}(E)$. X.W.~Chen \cite{chen} introduced a method
of constructing simple $L_{K}(E)$-modules from tail-equivalence classes of
infinite paths in $E$, and from the set of all finite paths ending at a sink.
This construction was extended in \cite{ar} and \cite{r-simple} to infinite
emitters and, by means of a twist, to cycles, producing several types of
simple modules, now usually called \emph{Chen modules}. In \cite{ar} it was
shown that, when $E$ is a finite graph in which every vertex is the base of at
most one cycle, these are the only simple modules up to isomorphism. An
alternative description of these modules was given in \cite{an}, and for the
rose graph $R_{n}$ with $n\geq2$ petals a number of further simple
$L_{K}(R_{n})$-modules were constructed by Nam, Srivastava and Vien \cite{nsv}.

In this paper we study simple $L_{K}(E)$-modules of a different kind, namely
modules generated by \emph{indicator functions} on sets of paths in $E$. The
initial step in this direction was taken by the first author in her
Ph.D.\ thesis \cite{pacheco}, where it was shown, using symbolic dynamics and
topological methods, that if $E$ is a finite graph and $S$ is a set of
infinite paths in $E$ with certain special properties, then the right
$L_{K}(E)$-module generated by the indicator function $\mathbf{1}_{S}$ is
simple. Here, we significantly improve the above result and its proof by
taking $E$ to be an arbitrary graph, allowing the sets $S$ to contain both
finite and infinite paths, and using purely algebraic arguments in our proofs.

\subsection*{The setting}

Let $D$ be the set of \emph{boundary paths} of $E$, that is, all infinite
paths together with all finite paths ending at a singular vertex. The space
$K^{D}$ of $K$-valued functions on $D$ is a right $L_{K}(E)$-module in a
natural way (Section~2). Recall \cite{chen} that two infinite paths
$p=e_{1}e_{2}e_{3}\cdots$ and $q=f_{1}f_{2}f_{3}\cdots$ are
\emph{tail-equivalent} if there are integers $m,n$ such that $e_{m+i}=f_{n+i}$
for all $i\geq1$. A tail-equivalence class $X=[p]$, and also the set $X$ of
all finite paths ending at a fixed sink, has the following property:

\begin{itemize}
\item[$(\#)$] If $\alpha$ is a finite path and $p$ is a path with
$r(\alpha)=s(p)$, then $\alpha p\in X$ if and only if $p\in X$.
\end{itemize}

\noindent

Thus $X$ is closed under adding and deleting initial segments. We call a
non-empty subset $S$ of $D$ \emph{admissible} if it satisfies $(\#)$ and is
strongly connected, that is, any two vertices lying on paths in $S$ are joined
by a path (Definition~\ref{def:adm}). For an admissible set $S$ we consider
the submodule
\[
N_{S}=\mathbf{1}_{S}\cdot L_{K}(E)\subseteq K^{D},
\]
together with the set $S^{0}$ of vertices lying on the paths in $S$ and the
hereditary saturated set $H=E^{0}\setminus S^{0}$ (Lemma~\ref{lem:basic}).
Throughout, one has to distinguish carefully between the \emph{paths} in $S$
and the \emph{vertices} in $S^{0}$ (Remark~\ref{important}).

\subsection*{Simplicity}

Our first main result, Theorem~\ref{thm:simple}, states that $N_{S}$ is a
simple module if and only if $S$ contains no breaking vertex of $H$;
equivalently, no finite path in $S$ ends at a breaking vertex of $H$. In
particular, $N_{S}$ is simple whenever $S$ consists only of infinite paths or
$E$ is row-finite (Corollaries~\ref{cor:infinite} and \ref{Finite case}), and
Example~\ref{ex:toeplitz} shows that the condition on breaking vertices cannot
be dropped. Up to isomorphism, a simple module $N_{S}$ depends only on the
vertex set $S^{0}$ and not on the paths in $S$ (Proposition~\ref{prop:indep}).
For instance, over the rose graph $R_{n}$, a single tail-equivalence class, a
union of several classes, and the set of all infinite paths all give the same
simple module (Example~\ref{ex:rose}).

\subsection*{Structure of $N_{S}$}

An admissible set $S$ is finite only in two degenerate situations: either $S$
is a single vertex, or $S$ is the set of rotations of $c^{\infty}$ for a cycle
$c$ that no other edge enters (Lemma~\ref{lem:finiteS}); in these cases
$N_{S}$ is finite dimensional. Apart from the case of a single vertex, the
module $N_{S}$ is far from being graded: $\operatorname{Hom}_{L_{K}(E)}
(N_{S},M)=0$ for every $\mathbb{Z}$-graded right $L_{K}(E)$-module $M$
(Theorem~\ref{thm:hom}). Consequently $N_{S}$ is not a graded module,
$\operatorname{Hom}_{L_{K}(E)}(N_{S},L_{K}(E))=0$, and $N_{S}$ is not
projective (Corollary~\ref{cor:hom}). If $S$ is infinite and $N_{S}$ is
simple, then the annihilator of $N_{S}$ is the graded primitive ideal
$I(H,B_{H})$ (Proposition~\ref{prop:ann}).

\subsection*{Comparison with Chen modules}

Since the Chen modules are left modules, in Section~\ref{sec:left} we let
$L_{K}(E)$ act on $K^{D}$ on the left through the involution of $L_{K}(E)$.
The resulting left module $L_{K}(E)\mathbf{1}_{S}$ has the same underlying
space and the same submodules as $N_{S}$, so it is simple exactly when $N_{S}$
is (Proposition~\ref{prop:leftstar}). When $S$ is infinite and contains no
breaking vertex of $H$, the simple left module $L_{K}(E)\mathbf{1}_{S}$ is not
isomorphic to any Chen simple module (Theorem~\ref{thm:chen}). For the rose
graph $R_{n}$ with one vertex $v$ and $n\geq2$ loops $e_{1},\dots,e_{n}$ we
obtain complete information, answering a question of T.G.~Nam. For every
admissible set $S$ we have $N_{S}\cong L_{K}(R_{n})/(v-\sum_{i}e_{i}^{\ast
})L_{K}(R_{n})$ (Proposition~\ref{prop:rose}). Moreover, the simple left
module $L_{K}(R_{n})\mathbf{1}_{S}$ is isomorphic to the twisted Chen module
$V_{[c^{\infty}]}^{P}$ of \cite{nsv}, where $c$ is a simple closed path and
$P\in GL_{n}(K)$, if and only if $c=e_{j}$ for some $j$ and the $j$-th column
of $P$ consists of $1$'s; and $L_{K}(R_{n})\mathbf{1}_{S}$ is not isomorphic
to any other simple module in the list of \cite[Theorem 4.7]{nsv}
(Theorem~\ref{thm:rose}).

\subsection*{Finite presentation}

If $S$ contains no breaking vertex of $H$ and is not a single vertex, then
$N_{S}$ is finitely presented if and only if $S^{0}$ is finite and only
finitely many edges of $E$ have both source and range in $S^{0}$; in that case
we give an explicit presentation of $N_{S}$ (Theorem~\ref{thm:fp}). The case
of a single vertex is settled in Proposition~\ref{prop:acyclic}. As a
consequence, over a finite graph every simple module $N_{S}$ is finitely
presented, whereas, by \cite{ar}, other simple modules over such a graph need
not be; and over a row-finite graph every simple module $N_{S}$ is finitely
presented if and only if $S^{0}$ is finite for every admissible set $S$ which
is not a single vertex (Corollary~\ref{cor:allfp}).

\subsection*{Organization of the paper}

Section~2 collects the preliminaries and introduces the module $K^{D}$.
Section~3 develops the basic properties of admissible sets. In Section~4 we
prove the simplicity criterion for $N_{S}$ and give examples. Section~5
studies homomorphisms from $N_{S}$ to graded modules, the annihilator of
$N_{S}$, and finite presentation of $N_{S}$. Section~\ref{sec:left} introduces
the left module $L_{K}(E)\mathbf{1}_{S}$, compares it with the Chen modules,
and treats the rose graph $R_{n}$ in detail. In Section~\ref{sec:cohn} we
consider the modules generated by indicator functions over the Cohn path
algebra $C_{K}(E)$, and determine when they are simple.

\section{Preliminaries}

For the general notation, terminology and results on Leavitt path algebras we
refer to \cite{aas}. We outline the concepts that we need.

A \emph{(directed) graph} $E=(E^{0},E^{1},r,s)$ consists of two sets $E^{0}$
and $E^{1}$ together with maps $r,s:E^{1}\rightarrow E^{0}$. The elements of
$E^{0}$ are called \emph{vertices} and the elements of $E^{1}$ \emph{edges}.
No restriction is placed on the number of vertices or on the number of edges
emitted by a vertex.

A vertex $v$ is a \emph{sink} if it emits no edges, a \emph{regular vertex} if
it emits a non-empty finite set of edges, and an \emph{infinite emitter} if it
emits infinitely many edges. A vertex which is a sink or an infinite emitter
is called \emph{singular}. A vertex is a \emph{source} if it is not the range
of any edge, and an \emph{isolated vertex} if it is both a source and a sink.
For each $e\in E^{1}$ we call $e^{\ast}$ a \emph{ghost edge}, and we let
$r(e^{\ast})=s(e)$ and $s(e^{\ast})=r(e)$. A \emph{path} $\mu$ of length
$|\mu|=n>0$ is a finite sequence of edges $\mu=e_{1}e_{2}\cdots e_{n}$ with
$r(e_{i})=s(e_{i+1})$ for $i=1,\dots,n-1$; we put $s(\mu)=s(e_{1})$,
$r(\mu)=r(e_{n})$ and $\mu^{\ast}=e_{n}^{\ast}\cdots e_{1}^{\ast}$. A vertex
$v$ is a path of length $0$ with $s(v)=r(v)=v$ and $v^{\ast}=v$. An
\emph{infinite path} is a sequence $p=e_{1}e_{2}e_{3}\cdots$ with
$r(e_{i})=s(e_{i+1})$ for all $i$; then $s(p)=s(e_{1})$. For a path $p$
(finite or infinite), $p^{0}$ denotes the set of vertices on $p$. If $\beta$
is a finite path and $p$ is a path, we write $\beta\preceq p$ if $\beta$ is an
initial segment of $p$, that is, $p=\beta p^{\prime}$ for some path
$p^{\prime}$ with $s(p^{\prime})=r(\beta)$; here $p^{\prime}$ may be the
vertex $r(\beta)$.

A finite path $\mu$ of positive length is \emph{closed} if $r(\mu)=s(\mu)$, in
which case $\mu$ is said to be \emph{based at} $s(\mu)$. A closed path
$\mu=e_{1}\cdots e_{n}$ is a \emph{cycle} if $s(e_{i})\neq s(e_{j})$ for
$i\neq j$. If $c$ is a closed path, $c^{\infty}$ denotes the infinite path
$ccc\cdots$. An \emph{exit} for a path $\mu=e_{1}\cdots e_{n}$ is an edge $f$
with $s(f)=s(e_{i})$ for some $i$ and $f\neq e_{i}$. The graph $E$ satisfies
\emph{Condition (L)} if every cycle in $E$ has an exit. If there is a path
from a vertex $u$ to a vertex $v$ we write $u\geq v$.

A subset $H$ of $E^{0}$ is \emph{hereditary} if $v\in H$ and $v\geq w$ imply
$w\in H$, and \emph{saturated} if, for every regular vertex $v$,
$r(s^{-1}(v))\subseteq H$ implies $v\in H$. For a hereditary saturated set
$H$, a \emph{breaking vertex} of $H$ is an infinite emitter $v\notin H$ such
that $s^{-1}(v)\cap r^{-1}(E^{0}\setminus H)$ is non-empty and finite. The set
of breaking vertices of $H$ is denoted by $B_{H}$.

\begin{definition}
\label{def:LPA}The \emph{Leavitt path algebra} $L_{K}(E)$ is the $K$-algebra
generated by a set $\{v:v\in E^{0}\}$ of pairwise orthogonal idempotents
together with a set of variables $\{e,e^{\ast}:e\in E^{1}\}$ satisfying:

\begin{enumerate}
\item[(1)] $s(e)e=e=er(e)$ for all $e\in E^{1}$;

\item[(2)] $r(e)e^{\ast}=e^{\ast}=e^{\ast}s(e)$ for all $e\in E^{1}$;

\item[(3)] (CK-1) $e^{\ast}e=r(e)$ and $e^{\ast}f=0$ if $e\neq f$, for all
$e,f\in E^{1}$;

\item[(4)] (CK-2) $v=\sum_{e\in s^{-1}(v)}ee^{\ast}$ for every regular vertex
$v\in E^{0}$.
\end{enumerate}
\end{definition}

Every element of $L_{K}(E)$ can be written as $\sum_{i=1}^{n}k_{i}\alpha
_{i}\beta_{i}^{\ast}$ where $k_{i}\in K$ and $\alpha_{i},\beta_{i}$ are finite
paths with $r(\alpha_{i})=r(\beta_{i})$. The algebra $L_{K}(E)$ is
$\mathbb{Z}$-graded, $L_{K}(E)=\bigoplus_{n\in\mathbb{Z}}L_{n}$, where $L_{n}$
is spanned by the monomials $\alpha\beta^{\ast}$ with $|\alpha|-|\beta|=n$. It
has a $K$-linear involution $a\mapsto a^{\ast}$ with $(\alpha\beta^{\ast
})^{\ast}=\beta\alpha^{\ast}$.

A right module $M$ over a ring $R$ is \emph{unitary} if $MR=M$. A right
$L_{K}(E)$-module $M$ is \emph{$\mathbb{Z}$-graded} if $M=\bigoplus
_{n\in\mathbb{Z}}M_{n}$ as abelian groups with $M_{n}L_{k}\subseteq M_{n+k}$
for all $n,k\in\mathbb{Z}$. An ideal $I$ of $L_{K}(E)$ is \emph{graded} if
$I=\bigoplus_{n}(I\cap L_{n})$. The annihilator of $M$ is the ideal
$\operatorname{ann}(M)=\{a\in L_{K}(E):Ma=0\}$.

For a hereditary saturated set $H$ and $v\in B_{H}$ put
\[
v^{H}=v-\sum_{e\in s^{-1}(v),\ r(e)\notin H}ee^{\ast},
\]
and let $I(H,B_{H})$ be the ideal of $L_{K}(E)$ generated by $H\cup
\{v^{H}:v\in B_{H}\}$. We shall use the following two known facts; see
\cite[Theorem 5.7 and Corollary 6.10]{T}.

\begin{theorem}
[\cite{T}]\label{thm:known}

\begin{enumerate}
\item[(a)] Let $H$ be a hereditary saturated subset of $E^{0}$ and let $F$ be
the graph with $F^{0}=E^{0}\setminus H$, $F^{1}=\{e\in E^{1}:r(e)\notin H\}$,
and $r,s$ restricted from $E$. Then $I(H,B_{H})$ is a graded ideal,
$I(H,B_{H})\cap E^{0}=H$, and there is an isomorphism of $K$-algebras from
$L_{K}(E)/I(H,B_{H})$ onto $L_{K}(F)$ sending $v+I(H,B_{H})$ to $v$ for each
$v\in E^{0}\setminus H$.

\item[(b)] If a graph $F$ satisfies Condition (L), then every non-zero ideal
of $L_{K}(F)$ contains a vertex.
\end{enumerate}
\end{theorem}

\subsection*{The boundary path module}

Let $D_{\infty}$ be the set of all infinite paths in $E$, let $D_{\mathrm{fin}%
}$ be the set of all finite paths $p$ (possibly of length $0$) such that
$r(p)$ is a singular vertex, and let
\[
D=D_{\infty}\cup D_{\mathrm{fin}} .
\]
Thus a vertex $v$ belongs to $D$ exactly when $v$ is singular. If $p\in D$ and
$\alpha$ is a finite path with $r(\alpha)=s(p)$, then $\alpha p\in D$; and if
$p\in D$ and $p=\alpha p^{\prime}$, then $p^{\prime}\in D$.

We shall work with right modules. The right $L_{K}(E)$-module $K^{D}$ of all
$K $-valued functions on $D$ is most easily obtained as the dual of a left
$L_{K}(E)$-module $KD$, which we describe first; the left module $KD$ is used
only for this purpose.

Let $KD$ be the $K$-vector space with basis $D$. For $v\in E^{0}$, $e\in
E^{1}$ and $p\in D$ define
\begin{align*}
v\cdot p  &  =%
\begin{cases}
p & \text{if }s(p)=v,\\
0 & \text{otherwise;}%
\end{cases}
\qquad e\cdot p=%
\begin{cases}
ep & \text{if }r(e)=s(p),\\
0 & \text{otherwise;}%
\end{cases}
\\[4pt]
e^{\ast}\cdot p  &  =%
\begin{cases}
p^{\prime} & \text{if }p=ep^{\prime},\\
0 & \text{otherwise.}%
\end{cases}
\end{align*}

\begin{lemma}
\label{lem:KD} The rules above make $KD$ a left $L_{K}(E)$-module. For finite
paths $\alpha,\beta$ with $r(\alpha)=r(\beta)$ and $p\in D$,
\begin{equation}
\label{eq:action}\alpha\beta^{\ast}\cdot p=%
\begin{cases}
\alpha p^{\prime} & \text{if }p=\beta p^{\prime},\\
0 & \text{if }\beta\not \preceq p.
\end{cases}
\end{equation}

\end{lemma}

\begin{proof}
Let $T_{v},T_{e},T_{e^{\ast}}$ be the linear operators on $KD$ defined by the
three rules. Clearly $T_{u}T_{v}=\delta_{u,v}T_{v}$, $T_{s(e)}T_{e}%
=T_{e}=T_{e}T_{r(e)}$ and $T_{r(e)}T_{e^{\ast}}=T_{e^{\ast}}=T_{e^{\ast}%
}T_{s(e)}$. For $e,f\in E^{1}$ and $p\in D$, $T_{e^{\ast}}T_{f}(p)$ is $0$
unless $r(f)=s(p)$ and $e=f$, in which case it equals $p$; hence $T_{e^{\ast}%
}T_{f}=\delta_{e,f}T_{r(e)}$. Let $v$ be a regular vertex and $p\in D$. If
$s(p)\neq v$, then $T_{v}(p)=0=\sum_{e\in s^{-1}(v)}T_{e}T_{e^{\ast}}(p)$. If
$s(p)=v$, then $p\neq v $ because $v$ is not singular, so $p=e_{0}p^{\prime}$
for a unique edge $e_{0}\in s^{-1}(v)$, and $\sum_{e\in s^{-1}(v)}%
T_{e}T_{e^{\ast}}(p)=e_{0}p^{\prime}=p=T_{v}(p)$. Thus the operators satisfy
the defining relations of $L_{K}(E)$ and, by the universal property of
$L_{K}(E)$, $KD$ is a left $L_{K}(E)$-module. Formula \eqref{eq:action}
follows by applying the generators one at a time.
\end{proof}

Note that in \eqref{eq:action} the value is $0$ in particular when $p$ is a
finite path which is a proper initial segment of $\beta$.

Let $K^{D}$ be the space of all functions $f:D\rightarrow K$; we identify it
with the dual space $\operatorname{Hom}_{K}(KD,K)$. Then $K^{D}$ is a right
$L_{K}(E)$-module under
\[
(f\cdot a)(p)=f(a\cdot p)\qquad(f\in K^{D},\ a\in L_{K}(E),\ p\in D),
\]
where $f$ is extended linearly to $KD$. Explicitly, by \eqref{eq:action},
\begin{equation}
(f\cdot\alpha\beta^{\ast})(p)=%
\begin{cases}
f(\alpha p^{\prime}) & \text{if }p=\beta p^{\prime},\\
0 & \text{if }\beta\not \preceq p.
\end{cases}
\label{eq:right}%
\end{equation}
The \emph{support} of $f\in K^{D}$ is the set $\operatorname{supp}(f)=\{p\in
D:f(p)\neq0\}$.

\textbf{Note}: In general $K^{D}$ is not a unitary $L_{K}(E)$-module. Indeed,
by \eqref{eq:right}, $f\cdot\alpha\beta^{\ast}$ vanishes outside the set of
paths starting at $s(\beta)$, so every element of $K^{D}\cdot L_{K}(E)$
vanishes outside the set of paths starting at finitely many vertices. Since
every vertex is the source of at least one boundary path, $\mathbf{1}%
_{D}\notin K^{D}\cdot L_{K}(E)$ when $E^{0}$ is infinite. (If $E^{0}$ is
finite, then $\sum_{v\in E^{0}}v$ is the identity of $L_{K}(E)$ and $K^{D}$ is unitary.)

For $X\subseteq D$, $\mathbf{1}_{X}\in K^{D}$ denotes the indicator function
of $X$, and for any $p\in D$, we use the notation $\delta_{p}$ to denote
$\mathbf{1}_{\{p\}}$.

\section{Admissible sets of paths}

\begin{definition}
\label{def:adm} Let $S$ be a subset of $D$.

\begin{enumerate}
\item[(i)] $S$ satisfies $(\#)$ if for every finite path $\alpha$ and every
$p\in D$ with $r(\alpha)=s(p)$ we have $\alpha p\in S$ if and only if $p\in S$.

\item[(ii)] For a finite path $\beta$ and a path $p$, $\beta\preceq p$ means
that $\beta$ is an initial segment of $p$, that is, $p=\beta p^{\prime}$ for
some path $p^{\prime}$ (possibly $p^{\prime}=r(\beta)$).

\item[(iii)] $S^{0}=\bigcup_{p\in S}p^{0}$ is the set of all vertices lying on
paths in $S$. For $v\in E^{0}$ and a finite path $\beta$ we put
\begin{gather*}
vS=\{p\in S:s(p)=v\},\qquad\beta S=\{p\in S:\beta\preceq p\},\\
s_{S}(v)=\{e\in E^{1}:s(e)=v,\ r(e)\in S^{0}\}.
\end{gather*}

\item[(iv)] $S$ is \emph{strongly connected} if $u\geq v$ for all $u,v\in
S^{0}$.

\item[(v)] $S$ is \emph{admissible} if it is non-empty, strongly connected and
satisfies $(\#)$.
\end{enumerate}
\end{definition}

\begin{remark}
\label{important} Two different things must be kept apart throughout: the
\emph{elements} of $S$, and the vertices in $S^{0} $. The set $S$ is a set of
paths, and \textquotedblleft$S$ contains $x$\textquotedblright\ always means
that the path $x$ is an element of $S$. A vertex $u$ is an element of $S$ only
if $u$, regarded as a path of length $0$, belongs to $S$; whereas $u\in S^{0}$
means only that $u$ lies on some path belonging to $S$.

(1) If $S$ consists of infinite paths, then no vertex and no edge is an
element of $S$, although $S^{0}$ is not empty. For example, if $E$ has one
vertex $v$ and two loops, and $S$ is the set of all infinite paths, then
$v\notin S$ while $S^{0}=\{v\}$.

(2) Let $S$ satisfy $(\#)$ and let $p=e_{1}e_{2}\cdots e_{n}$ be a finite path
belonging to $S$. Deleting initial segments, we see that the terminal
segments
\[
e_{1}e_{2}\cdots e_{n},\quad e_{2}\cdots e_{n},\quad\dots,\quad e_{n},\quad
r(p)
\]
of $p$ are elements of $S$; in particular the vertex $r(p)$ is an element of
$S$. It does not follow that each edge $e_{k}$ is an element of $S$: for $k<n
$, $e_{k}\in S$ holds if and only if the vertex $r(e_{k})$ is an element of
$S$. All the vertices on $p$ belong to $S^{0}$.

(3) Consequently, for a set $S$ satisfying $(\#)$, a vertex $u$ is an element
of $S$ if and only if some finite path in $S$ ends at $u$; and $S$ contains no
vertex if and only if $S\subseteq D_{\infty}$.

(4) \textquotedblleft$S$ is finite\textquotedblright\ means that $S$ has only
finitely many elements, that is, $S$ consists of finitely many paths. This is
not the same as $S^{0}$ being finite: in the example in (1), $S^{0}$ has one
element and $S$ is infinite.
\end{remark}

\begin{lemma}
\label{lem:basic} Let $S\subseteq D$ be non-empty and satisfy $(\#)$.

\begin{enumerate}
\item[(a)] If $p\in S$ and $p=\alpha p^{\prime}$ with $\alpha$ a finite path,
then $p^{\prime}\in S$. In particular, if $p\in S$ is a finite path then the
vertex $r(p)$ belongs to $S$, and every finite path ending at $r(p)$ belongs
to $S$. Moreover $S\cap D_{\infty}$ is a union of tail-equivalence classes.

\item[(b)] A vertex $v$ belongs to $S^{0}$ if and only if $vS\neq\emptyset$.

\item[(c)] If $u\in E^{0}$ and $u\geq v$ for some $v\in S^{0}$, then $u\in
S^{0}$. Hence $H=E^{0}\setminus S^{0}$ is hereditary.

\item[(d)] $H=E^{0}\setminus S^{0}$ is saturated.

\item[(e)] If $r(\beta)\in S^{0}$ for a finite path $\beta$, then $\beta
S\neq\emptyset$; if $r(\beta)\notin S^{0}$, then $\beta S=\emptyset$.
\end{enumerate}
\end{lemma}

\begin{proof}
(a) The first statement is $(\#)$ applied to $\alpha$ and $p^{\prime}$. If
$p\in S$ is finite, write $p=p\,r(p)$ to get $r(p)\in S$, and then
$\alpha=\alpha\, r(p)\in S$ for every finite path $\alpha$ ending at $r(p)$.
If $p\in S\cap D_{\infty}$ and $q$ is tail-equivalent to $p$, then $p=\alpha
x$ and $q=\beta x$ for some infinite path $x$ and finite paths $\alpha,\beta$;
so $x\in S$ and then $q\in S$.

(b) If $v\in S^{0}$, then $v$ lies on some $p\in S$, so $p=\alpha p^{\prime} $
with $s(p^{\prime})=v$ (if $p$ is finite and $v=r(p)$, take $p^{\prime}=v $).
By (a), $p^{\prime}\in vS$. The converse is clear.

(c) Let $p^{\prime}\in vS$ and let $\mu$ be a path from $u$ to $v$. Then $\mu
p^{\prime}\in S$ by $(\#)$, so $u\in S^{0}$.

(d) Let $v$ be a regular vertex with $r(s^{-1}(v))\subseteq H$ and suppose
$v\in S^{0}$. Choose $p\in vS$. Since $v$ is regular, $p\neq v$, so $p$ begins
with an edge $e\in s^{-1}(v)$, and $r(e)\in p^{0}\subseteq S^{0}$, a contradiction.

(e) If $r(\beta)\in S^{0}$, choose $p^{\prime}\in r(\beta)S$; then $\beta
p^{\prime}\in\beta S$. If $p=\beta p^{\prime}\in S$, then $r(\beta)\in
p^{0}\subseteq S^{0}$.
\end{proof}

\begin{lemma}
\label{lem:single} Let $S$ be admissible and suppose that $S$ contains a
vertex $w$ with $s_{S}(w)=\emptyset$ (for instance, a sink). Then $S=\{w\}$.
\end{lemma}

\begin{proof}
Let $u\in S^{0}$. Then $w\geq u$. If $u\neq w$, a path from $w$ to $u$ begins
with an edge $e$ and $r(e)\geq u$, so $r(e)\in S^{0}$ by Lemma~\ref{lem:basic}%
(c), contradicting $s_{S}(w)=\emptyset$. Hence $S^{0}=\{w\}$. A path of
positive length all of whose vertices equal $w$ begins with a loop at $w$,
which would belong to $s_{S}(w)$. Hence $S=\{w\}$.
\end{proof}

Note that $S=\{w\}$ is admissible exactly when $w$ is a singular vertex which
is a source; that is, when $w$ is an isolated vertex, or $w$ is an infinite
emitter which is a source.

We next describe the admissible sets that are finite.

\begin{lemma}
\label{lem:finiteS} For an admissible set $S$ the following are equivalent:

\begin{enumerate}
\item[(i)] $S$ is a finite set, that is, $S$ consists of only finitely many
paths (a vertex belonging to $S$ being counted as a path of length $0$);

\item[(ii)] $|vS|=1$ for some $v\in S^{0}$;

\item[(iii)] $|vS|=1$ for every $v\in S^{0}$;

\item[(iv)] either $S=\{w\}$ for a vertex $w$ (which is then an isolated
vertex or an infinite emitter which is a source), or there is a cycle $c$ in
$E$ such that the only edges of $E$ with range in $c^{0}$ are the edges of $c
$, and $S$ consists of the $|c|$ infinite paths $c_{u}^{\infty}$, $u\in c^{0}
$, where $c_{u}$ denotes the rotation of $c$ based at $u$.
\end{enumerate}

If these conditions fail, then $vS$ is infinite for every $v\in S^{0}$.
\end{lemma}

\begin{proof}
Recall that $vS\neq\emptyset$ for $v\in S^{0}$. If $u,v\in S^{0}$ and
$\lambda$ is a path from $u$ to $v$, then $y\mapsto\lambda y$ is an injective
map from $vS$ to $uS$. Hence $|vS|\leq|uS|$ for all $u,v\in S^{0}$, so all the
sets $vS$, $v\in S^{0}$, have the same cardinality if one of them is finite.
This proves (ii)$\Leftrightarrow$(iii). Clearly (iv)$\Rightarrow$(i).

(iii)$\Rightarrow$(iv). Suppose first that $S$ contains a finite path. Then
$S$ contains a vertex $w$ by Lemma~\ref{lem:basic}(a), and $wS=\{w\}$. If $e$
is an edge with $r(e)=w$, then $e\in S$ by $(\#)$, so $s(e)\in S^{0}$ and
there is a path $\mu$ from $w$ to $s(e)$; but then $\mu e\in wS$ has positive
length, a contradiction. So no edge has range $w$. As every $u\in S^{0}$
satisfies $u\geq w$, we get $S^{0}=\{w\}$, and then $S=wS=\{w\}$.

Suppose now that $S\subseteq D_{\infty}$. Fix $p\in S$ and let $v=s(p)$, so
$vS=\{p\}$. The second vertex of $p$ lies in $S^{0}$, so by strong
connectivity there are closed paths based at $v$. Let $c$ be one of minimal
length. Then $c$ is a cycle, since otherwise deleting a closed subpath of $c$
would give a shorter closed path based at $v$. For every closed path $d$ based
at $v$ we have $dp\in vS$, so $dp=p$ and hence $p=d^{\infty}$. In particular
$p=c^{\infty}$. Let $u\in S^{0}$, $u\neq v$, and choose paths $\lambda$ from
$v$ to $u$ and $\mu$ from $u$ to $v$. Then $(\lambda\mu)^{\infty}=c^{\infty}$,
so $\lambda\preceq c^{\infty}$ and $u=r(\lambda)\in c^{0}$. Thus $S^{0}=c^{0}%
$, and for $u\in c^{0}$ the unique element of $uS$ is the tail $c_{u}^{\infty
}$ of $c^{\infty}$ (Lemma~\ref{lem:basic}(a)). Hence $S=\{c_{u}^{\infty}:u\in
c^{0}\}$. Finally, let $f$ be an edge with $r(f)\in c^{0}$. Then $s(f)\in
S^{0}=c^{0}$ by Lemma~\ref{lem:basic}(c), and $fc_{r(f)}^{\infty}\in
s(f)S=\{c_{s(f)}^{\infty}\}$, so $f$ is the first edge of $c_{s(f)}$, an edge
of $c$.

(i)$\Rightarrow$(ii), and the last statement. Suppose $|vS|\geq2$ for every
$v\in S^{0}$. Then $S$ contains a path of positive length: otherwise $S$ would
consist of vertices and $wS=\{w\}$ for $w\in S$. Let $e$ be an edge on a path
in $S$ and $v=s(e)$. Since $r(e)\geq v$, there is a closed path $d$ based at
$v$. Choose $y\in vS$ with $y\neq d^{\infty}$. The paths $d^{n}y$, $n\geq0$,
belong to $vS$ and are pairwise distinct: this is clear if $y$ is finite, and
if $y$ is infinite and $d^{n}y=d^{m}y$ with $n<m$, then $y=d^{m-n}y$, which
forces $y=d^{\infty}$. So $vS$ is infinite, and by the first paragraph $uS$ is
infinite for every $u\in S^{0}$.
\end{proof}

Thus an admissible set $S$ is either a finite set of paths, in which case it
is one of the sets described in Lemma~\ref{lem:finiteS}(iv) and $S^{0}$ is
finite as well, or an infinite set of paths, in which case $vS$ is infinite
for every $v\in S^{0}$. Throughout, ``$S$ is finite'' and ``$S$ is infinite''
refer to the number of paths in $S$.

For an admissible set $S$ which is not a single vertex, the vertex set $S^{0}
$ is already determined by any closed path in it, as the next lemma shows.

\begin{lemma}
\label{lem:source} Let $S$ be an admissible set which is not of the form
$\{w\}$ with $w$ a vertex.

\begin{enumerate}
\item[(a)] $S$ contains a path of positive length, and $S^{0}$ contains a
vertex which is the base of a closed path all of whose vertices lie in $S^{0}
$.

\item[(b)] $S^{0}=\{u\in E^{0}:u\geq v\}$ for every $v\in S^{0}$, and every
vertex of $S^{0}$ lies on a closed path all of whose vertices lie in $S^{0}$.
In particular, every vertex of $S^{0}$ is the source of at least one edge
whose range is in $S^{0}$.

\item[(c)] If $d$ is a closed path based at a vertex $v\in S^{0}$, then the
tail-equivalence class $[d^{\infty}]$ is an admissible set of infinite paths
with $[d^{\infty}]^{0}=S^{0}$.
\end{enumerate}
\end{lemma}

\begin{proof}
(a) If $S$ consisted of vertices only, it would contain two distinct vertices
$u,w$, and a path $\mu$ from $u$ to $w$ would satisfy $\mu=\mu w\in S$ by
$(\#)$, which is impossible. So $S$ contains a path of positive length. If $e$
is an edge on such a path, then $s(e),r(e)\in S^{0}$ and a path $\mu$ from
$r(e)$ to $s(e)$ gives a closed path $e\mu$ based at $s(e)$; all its vertices
lie in $S^{0}$ by Lemma~\ref{lem:basic}(c).

(b) If $v\in S^{0}$, then $\{u:u\geq v\}\subseteq S^{0}$ by
Lemma~\ref{lem:basic}(c), and $S^{0}\subseteq\{u:u\geq v\}$ by strong
connectivity. Let $c$ be a closed path in $S^{0}$ based at $w$, as in (a), and
let $u\in S^{0}$. Paths $\lambda$ from $u$ to $w$ and $\mu$ from $w$ to $u$
have all their vertices in $S^{0}$, and $\lambda c\mu$ is a closed path based
at $u$.

(c) A tail-equivalence class satisfies $(\#)$ and consists of infinite paths.
A vertex $u$ lies on a path tail-equivalent to $d^{\infty}$ if and only if
$u\geq v$: such a path eventually runs along $d$, and conversely, if $\lambda$
is a path from $u$ to $v$, then $\lambda d^{\infty}\in\lbrack d^{\infty}]$. By
(b), $[d^{\infty}]^{0}=\{u:u\geq v\}=S^{0}$, which is strongly connected; so
$[d^{\infty}]$ is admissible.
\end{proof}

\section{The module generated by an indicator function}

Recall from Definition~\ref{def:adm} that, for a finite path $\beta$ and a
path $p$, $\beta\preceq p$ means that $\beta$ is an initial segment of $p$,
and that $s_{S}(v)=\{e\in E^{1}:s(e)=v,\ r(e)\in S^{0}\}$ for a vertex $v$.

Throughout this section $S$ is a non-empty strongly connected subset of $D$
satisfying $(\#)$, and
\[
N_{S}=\mathbf{1}_{S}\cdot L_{K}(E)\subseteq K^{D}.
\]

\begin{lemma}
\label{lem:ind} Let $\alpha,\beta$ be finite paths with $r(\alpha)=r(\beta)$.
Then $\mathbf{1}_{S}\cdot\alpha\beta^{\ast}=\mathbf{1}_{\beta S}$. Consequently:

\begin{enumerate}
\item[(a)] $N_{S}=\operatorname{span}_{K}\{\mathbf{1}_{\beta S}:\beta\text{ a
finite path with }r(\beta)\in S^{0}\}$, every element of $N_{S}$ has support
contained in $S$, and $N_{S}$ is a unitary $L_{K}(E)$-module;

\item[(b)] $\mathbf{1}_{S}\cdot v=\mathbf{1}_{vS}$ for $v\in E^{0}$, and this
is non-zero exactly when $v\in S^{0}$;

\item[(c)] $\mathbf{1}_{vS}\cdot\lambda=\mathbf{1}_{r(\lambda)S}$ for every
finite path $\lambda$ with $s(\lambda)=v$;

\item[(d)] $\mathbf{1}_{r(\beta)S}\cdot\beta^{\ast}=\mathbf{1}_{\beta S}$ for
every finite path $\beta$.
\end{enumerate}
\end{lemma}

\begin{proof}
Let $p\in D$. By \eqref{eq:right}, $(\mathbf{1}_{S}\cdot\alpha\beta^{\ast
})(p)=1$ if and only if $p=\beta p^{\prime}$ and $\alpha p^{\prime}\in S$. By
$(\#)$, $\alpha p^{\prime}\in S\Leftrightarrow p^{\prime}\in S\Leftrightarrow
\beta p^{\prime}\in S$. So $(\mathbf{1}_{S}\cdot\alpha\beta^{\ast})(p)=1$ if
and only if $p\in\beta S$.

(a) The first two statements follow because the monomials $\alpha\beta^{\ast}
$ span $L_{K}(E)$ and $\beta S=\emptyset$ when $r(\beta)\notin S^{0}$
(Lemma~\ref{lem:basic}(e)). By \eqref{eq:right}, $\mathbf{1}_{\beta S}\cdot
s(\beta)=\mathbf{1}_{\beta S}$, so $N_{S}$ is unitary. (b) is the case
$\alpha=\beta=v$, together with Lemma~\ref{lem:basic}(b). (c) $\mathbf{1}%
_{vS}\cdot\lambda=\mathbf{1}_{S}\cdot v\lambda=\mathbf{1}_{S}\cdot
\lambda\,r(\lambda)^{\ast}=\mathbf{1}_{r(\lambda)S}$. (d) $\mathbf{1}%
_{r(\beta)S}\cdot\beta^{\ast}=\mathbf{1}_{S}\cdot r(\beta)\beta^{\ast
}=\mathbf{1}_{\beta S}$.
\end{proof}

\begin{remark}
Every element of $N_{S}$ is supported on finitely many of the sets $vS$. Hence
$\mathbf{1}_{S}$ itself belongs to $N_{S}$ if and only if $S^{0}$ is finite.
In general $N_{S}$ is the submodule of $K^{D}$ generated by the functions
$\mathbf{1}_{vS}$, $v\in S^{0}$.
\end{remark}

The next lemma expresses that the elements of $N_{S}$ are locally constant.

\begin{lemma}
\label{lem:local} Let $g=\sum_{i=1}^{t}k_{i}\mathbf{1}_{\beta_{i}S}\in N_{S} $
with $k_{i}\in K$, and let $\gamma$ be a finite path which is not a proper
initial segment of any $\beta_{i}$. Put $c_{\gamma}=\sum\{k_{i}:\beta
_{i}\preceq\gamma\}$. Then
\[
g\cdot\gamma=c_{\gamma}\,\mathbf{1}_{r(\gamma)S},\qquad\text{that is,}\qquad
g(\gamma y)=c_{\gamma}\ \text{ for all }y\in r(\gamma)S.
\]

\end{lemma}

\begin{proof}
Let $y\in D$. If $s(y)\neq r(\gamma)$ then $(g\cdot\gamma)(y)=0$. Let
$s(y)=r(\gamma)$, so $(g\cdot\gamma)(y)=g(\gamma y)$. If $y\notin S$ then
$\gamma y\notin S$ and $g(\gamma y)=0$. If $y\in S$ then $\gamma y\in S$ and
$g(\gamma y)=\sum\{k_{i}:\beta_{i}\preceq\gamma y\}$. Now $\beta_{i}$ and
$\gamma$ are both initial segments of $\gamma y$ whenever $\beta_{i}%
\preceq\gamma y$, so one of them is an initial segment of the other; as
$\gamma$ is not a proper initial segment of $\beta_{i}$, we get $\beta
_{i}\preceq\gamma y\Leftrightarrow\beta_{i}\preceq\gamma$. Hence $g(\gamma
y)=c_{\gamma}$.
\end{proof}

\begin{lemma}
\label{lem:key} Let $S$ be admissible, $S\neq\{w\}$ for every vertex $w$, and
suppose that $S$ contains no breaking vertex of $H=E^{0}\setminus S^{0}$. Then
for every non-zero $g\in N_{S}$ there exist a finite path $\gamma$ with
$r(\gamma)\in S^{0}$ and a non-zero $c\in K$ such that $g\cdot\gamma
=c\,\mathbf{1}_{r(\gamma)S}$; in particular $g(\gamma y)=c$ for all $y\in
r(\gamma)S$.
\end{lemma}

\begin{proof}
Write $g=\sum_{i=1}^{t}k_{i} \mathbf{1}_{\beta_{i}S}$ and let $m=\max
_{i}|\beta_{i}|$. Choose $p\in D$ with $g(p)\neq0$. Then $p\in S$ and
$g(p)=\sum\{k_{i}:\beta_{i}\preceq p\}$.

\emph{Case 1: $|p|\geq m$} (this includes the case that $p$ is infinite). Let
$\gamma$ be the initial segment of $p$ of length $m$. Then $\gamma$ is not a
proper initial segment of any $\beta_{i}$, and $\beta_{i}\preceq
\gamma\Leftrightarrow\beta_{i}\preceq p$ because $|\beta_{i}|\leq m$. So
$c_{\gamma}=g(p)\neq0$, $r(\gamma)\in p^{0}\subseteq S^{0}$, and
Lemma~\ref{lem:local} gives $g\cdot\gamma=g(p)\mathbf{1}_{r(\gamma)S}$.

\emph{Case 2: $|p|<m$.} Then $p$ is finite and $w=r(p)$ is a singular vertex
belonging to $S$ (Lemma~\ref{lem:basic}(a)). Since $S\neq\{w\}$,
Lemma~\ref{lem:single} shows that $s_{S}(w)\neq\emptyset$; so $w$ is an
infinite emitter, and since $w\notin B_{H}$ the set $s_{S}(w)$ is infinite.
For each $i$ there is at most one edge $e$ with $pe\preceq\beta_{i}$. Hence we
can choose $e\in s_{S}(w)$ such that $pe\not \preceq \beta_{i}$ for all $i$.
Put $\gamma=pe$. Then $\gamma$ is not an initial segment of any $\beta_{i}$,
so $\beta_{i}\preceq\gamma\Leftrightarrow\beta_{i}\preceq p$ and $c_{\gamma
}=g(p)\neq0 $. By Lemma~\ref{lem:local}, $g\cdot\gamma=g(p)\mathbf{1}_{r(e)S}%
$, and $r(e)\in S^{0}$.
\end{proof}

\begin{theorem}
\label{thm:simple} Let $E$ be an arbitrary graph, let $S\subseteq D$ be an
admissible set of paths and let $H=E^{0}\setminus S^{0}$. Then $N_{S}%
=\mathbf{1}_{S}\cdot L_{K}(E)$ is a simple right $L_{K}(E)$-module if and only
if $S$ contains no breaking vertex of $H$.
\end{theorem}

\begin{proof}
Note that $N_{S}\neq0$. If $S=\{w\}$ for a vertex $w$, then $s_{S}%
(w)=\emptyset$, so $w\notin B_{H}$, and $N_{S}=K\delta_{w}$ is one-dimensional
and therefore simple. So assume $S\neq\{w\}$ for every vertex $w$.

Suppose $S\cap B_{H}=\emptyset$ and let $0\neq g\in N_{S}$. By
Lemma~\ref{lem:key}, $\mathbf{1}_{v_{0}S}\in g\cdot L_{K}(E)$ for some
$v_{0}\in S^{0}$. Let $v\in S^{0}$. By strong connectivity there is a path
$\lambda$ from $v_{0}$ to $v$, and $\mathbf{1}_{vS}=\mathbf{1}_{v_{0}S}%
\cdot\lambda\in g\cdot L_{K}(E)$ by Lemma~\ref{lem:ind}(c). Then, for every
finite path $\beta$ with $r(\beta)\in S^{0}$, $\mathbf{1}_{\beta S}%
=\mathbf{1}_{r(\beta)S}\cdot\beta^{\ast}\in g\cdot L_{K}(E)$ by
Lemma~\ref{lem:ind}(d). By Lemma~\ref{lem:ind}(a), $g\cdot L_{K}(E)=N_{S}$.
Hence $N_{S}$ is simple.

Conversely, suppose that $S$ contains a breaking vertex $w$ of $H$, so $w\in
S\cap B_{H}$. Let $s_{S}(w)=\{e_{1},\dots,e_{k}\}$, $k\geq1$. The set $wS$ is
the union of $\{w\}$ and the sets $e_{1}S,\dots,e_{k}S$, because $w\in wS $
and every path $p\in wS$ with $p\neq w$ begins with an edge whose range lies
in $S^{0}$, that is, with one of the edges $e_{1},\dots,e_{k}$. As these sets
are pairwise disjoint, $\mathbf{1}_{wS}=\delta_{w}+\sum_{i=1}^{k}%
\mathbf{1}_{e_{i}S}$. Since $\mathbf{1}_{S}\cdot w=\mathbf{1}_{wS}$ and
$\mathbf{1}_{S}\cdot e_{i}e_{i}^{\ast}=\mathbf{1}_{e_{i}S}$ by
Lemma~\ref{lem:ind}, we get
\[
\delta_{w}=\mathbf{1}_{wS}-\sum_{i=1}^{k}\mathbf{1}_{e_{i}S}=\mathbf{1}
_{S}\cdot\Bigl(w-\sum_{i=1}^{k}e_{i}e_{i}^{\ast}\Bigr),
\]
so $\delta_{w}$ is a non-zero element of $N_{S}$. We claim that $\delta
_{w}\cdot L_{K}(E)$ is a proper submodule of $N_{S}$. By \eqref{eq:right}, for
finite paths $\alpha,\beta$ with $r(\alpha)=r(\beta)$, $(\delta_{w}\cdot
\alpha\beta^{\ast})(p)=1$ exactly when $p=\beta p^{\prime}$ and $\alpha
p^{\prime}=w$, that is, when $\alpha=w$, $p^{\prime}=w$ and $p=\beta$. So
$\delta_{w}\cdot\alpha\beta^{\ast}$ is $\delta_{\beta}$ if $\alpha=w$, and $0$
otherwise. Therefore
\[
\delta_{w}\cdot L_{K}(E)=\operatorname{span}_{K}\{\delta_{\beta}:\beta\text{ a
finite path with }r(\beta)=w\}
\]
consists of functions with finite support. On the other hand, $r(e_{1})\in
S^{0}$ gives a path $\mu$ from $r(e_{1})$ to $w$, and the closed path
$d=e_{1}\mu$ satisfies $d^{n}\in wS$ for all $n\geq0$ (Lemma~\ref{lem:basic}%
(a)). So $\mathbf{1}_{wS}\in N_{S}$ has infinite support and does not belong
to $\delta_{w}\cdot L_{K}(E)$. Thus $\delta_{w}\cdot L_{K}(E)$ is a non-zero
proper submodule of $N_{S}$, and $N_{S}$ is not simple.
\end{proof}

Note that the condition in Theorem~\ref{thm:simple} concerns the set of paths
$S$, not the set of vertices $S^{0}$. A breaking vertex $u$ of $H$ is by
definition a vertex not in $H$, so it always belongs to $S^{0}=E^{0}\setminus
H$. To say that $S$ contains no breaking vertex of $H$ means that no such $u$,
regarded as a path of length $0$, is an element of $S $; by
Lemma~\ref{lem:basic}(a) this is the same as saying that no finite path in $S$
ends at a breaking vertex of $H$.

Suppose now that $S\subseteq D_{\infty}$. The elements of $B_{H}$ are
vertices, that is, paths of length $0$, and $S$ contains no finite paths. So
$S\cap B_{H}=\emptyset$, and Theorem~\ref{thm:simple} gives:

\begin{corollary}
\label{cor:infinite} If $S\subseteq D_{\infty}$ is a non-empty strongly
connected set of infinite paths satisfying $(\#)$, then $N_{S}=\mathbf{1}%
_{S}\cdot L_{K}(E)$ is a simple $L_{K}(E)$-module.
\end{corollary}

Since a row-finite graph contains no infinite emitters, and, in particular, no
breaking vertices, we get the following corollary.

\begin{corollary}
\label{Finite case} If $E$ is a row-finite graph, then $N_{S}$ is a simple
module for any choice of admissible set $S\subseteq D$.
\end{corollary}

\begin{corollary}
\label{cor:support} Let $S$ be an infinite admissible set containing no
breaking vertex of $H$. Then every non-zero element of $N_{S}$ has infinite support.
\end{corollary}

\begin{proof}
By Lemma~\ref{lem:key}, a non-zero $g\in N_{S}$ is non-zero on $\{\gamma
y:y\in r(\gamma)S\}$ for some finite path $\gamma$ with $r(\gamma)\in S^{0}$,
and $r(\gamma)S$ is infinite by Lemma~\ref{lem:finiteS}.
\end{proof}

The simple module $N_{S}$ depends only on the set of vertices $S^{0}$ and not
on the paths connecting these vertices:

\begin{proposition}
\label{prop:indep} Let $S$ and $S^{\prime}$ be admissible sets with
$S^{0}=S^{\prime0}$, neither of which contains a breaking vertex of
$H=E^{0}\setminus S^{0}$. Then $N_{S}\cong N_{S^{\prime}}$.
\end{proposition}

\begin{proof}
The set $T=S\cup S^{\prime}$ satisfies $(\#)$, $T^{0}=S^{0}$ is strongly
connected, and $T\cap B_{H}=\emptyset$. So $N_{T}$ is simple by
Theorem~\ref{thm:simple}, and it suffices to show $N_{T}\cong N_{S}$. By
$(\#)$ and Lemma~\ref{lem:basic}(a), the subspace $KS$ of $KD$ is an
$L_{K}(E)$-submodule. Hence the restriction map $\rho:K^{D}\rightarrow K^{S}$,
$f\mapsto f|_{S}$, is a homomorphism of right $L_{K}(E)$-modules, where
$K^{S}=\operatorname{Hom}_{K}(KS,K)$. It is injective on functions supported
in $S$, and $\rho(\mathbf{1}_{T}\cdot a)=\rho(\mathbf{1}_{T})\cdot
a=\rho(\mathbf{1}_{S})\cdot a=\rho(\mathbf{1}_{S}\cdot a)$ for $a\in L_{K}%
(E)$. So $\mathbf{1}_{T}\cdot a\mapsto\mathbf{1}_{S}\cdot a$ is a well-defined
surjective homomorphism $N_{T}\rightarrow N_{S}$. It is non-zero, and $N_{T}$
is simple, so it is an isomorphism from $N_{T}$ to $N_{S}$. Similarly,
$N_{T}\cong N_{S^{\prime}}$.
\end{proof}

\begin{example}
\label{ex:rose} Let $E=R_{n}$ be the graph with one vertex $v$ and $n\geq2$
loops $e_{1},\dots,e_{n}$. Here $D=D_{\infty}$. Every non-empty union $S$ of
tail-equivalence classes of infinite paths is admissible with $S^{0}=\{v\} $,
and it is infinite. By Corollary~\ref{cor:infinite} and
Proposition~\ref{prop:indep}, the modules $N_{S}$ are simple and isomorphic to
one another. For instance $S$ may be a single class $[p]$, the union
$[p]\cup\lbrack q]$ of two different classes, or all of $D_{\infty}$; in the
last case $N_{S}$ is the space of locally constant functions on $D_{\infty}$.
Thus an admissible set need not be a single tail-equivalence class.
\end{example}

\begin{example}
\label{ex:emitter} Let $E$ be the graph with a vertex $w$, loops $e_{2i}$
($i\in\mathbb{N}$) based at $w$, sinks $z_{i}$ ($i\in\mathbb{N}$), and edges
$e_{2i+1}$ from $w$ to $z_{i}$. Thus $w$ is an infinite emitter. Let $S_{1}$
be the set of all finite paths ending at $w$ (the vertex $w$ and all finite
words in the loops $e_{2i}$), let $S_{2}$ be the set of all infinite paths in
$E$ (the infinite words in the loops), and let $S_{3}=S_{1}\cup S_{2}$. Each
$S_{j}$ is admissible and infinite, with $S_{j}^{0}=\{w\}$ and $H=\{z_{i}%
:i\in\mathbb{N}\}$. Since $s_{S_{j}}(w)=\{e_{2i}:i\in\mathbb{N}\}$ is
infinite, $B_{H}=\emptyset$. By Theorem~\ref{thm:simple} and
Proposition~\ref{prop:indep}, $N_{S_{1}}\cong N_{S_{2}}\cong N_{S_{3}}$ is a
simple module. Note that the smaller set $\{w\}\cup\{e_{2i}:i\in\mathbb{N}\}$
does not satisfy $(\#)$, since it contains $e_{4}$ but not $e_{2}e_{4}$.
\end{example}

\begin{example}
\label{ex:toeplitz} The condition on breaking vertices in
Theorem~\ref{thm:simple} cannot be omitted. Let $E$ be the graph with a vertex
$w$, a single loop $e$ based at $w$, sinks $z_{i}$ ($i\in\mathbb{N}$) and
edges $f_{i}$ from $w$ to $z_{i}$. Then $w$ is an infinite emitter. The set
$S=\{e^{n}:n\geq0\}$ of all finite paths ending at $w$ (with $e^{0}=w$) is
admissible, and every path in $S$ ends at the infinite emitter $w$. Here
$H=\{z_{i}:i\in\mathbb{N}\}$, $s_{S}(w)=\{e\}$ and $w\in S\cap B_{H}$.
Identifying $e^{n}$ with $n$, $N_{S}$ is the space spanned by the functions
$\mathbf{1}_{[k,\infty)}$ on $\mathbb{N}$, that is, the space of eventually
constant functions, and $\delta_{w}\cdot L_{K}(E)$ is the proper submodule of
finitely supported functions. So $N_{S}$ is not simple.
\end{example}

\begin{example}
\label{ex:loop} Let $E$ be the graph with one vertex $v$ and one loop $c$, so
that $L_{K}(E)\cong K[x,x^{-1}]$. The only admissible set is $S=\{c^{\infty
}\}$, which is finite, and $N_{S}=K\mathbf{1}_{S}$ with $\mathbf{1}_{S}\cdot
c=\mathbf{1}_{S}=\mathbf{1}_{S}\cdot c^{\ast}$. So $N_{S}\cong K[x,x^{-1}%
]/(x-1)$. More generally, if $S=\{c_{u}^{\infty}:u\in c^{0}\}$ is as in
Lemma~\ref{lem:finiteS}(iv), then $N_{S}$ has the basis $\{\mathbf{1}%
_{uS}:u\in c^{0}\}$; for each edge $e$ of $c$ we have $\mathbf{1}_{s(e)S}\cdot
e=\mathbf{1}_{r(e)S}$ and $\mathbf{1}_{r(e)S}\cdot e^{\ast}=\mathbf{1}%
_{s(e)S}$, and every other product of a basis element with an edge or a ghost
edge is zero.
\end{example}

\begin{remark}
\label{rem:acyclic} If $S$ is admissible and contains a path of positive
length, then $E$ contains a cycle: an edge $e$ on a path in $S$ has
$s(e),r(e)\in S^{0}$, and $r(e)\geq s(e)$. Hence in an acyclic graph the only
admissible sets are the sets $\{w\}$, where $w$ is an isolated vertex or an
infinite emitter which is a source.
\end{remark}

\section{Properties of the module $N_{S}$}

Throughout this section, $N_{S}=\mathbf{1}_{S}\cdot L_{K}(E)$ is a simple
module with $S\subseteq D$ satisfying the conditions of Theorem
\ref{thm:simple} and $H=E^{0}\setminus S^{0}$.

\subsection*{Homomorphisms to graded modules}

\begin{theorem}
\label{thm:hom} Suppose $S$ is not of the form $\{w\}$, with $w$ a vertex.
Then $\operatorname{Hom}_{L_{K}(E)}(N_{S},M)=0$ for every $\mathbb{Z}$-graded
right $L_{K}(E)$-module $M$.
\end{theorem}

\begin{proof}
First, $S$ contains a path of positive length: otherwise $S$ would contain two
distinct vertices $u,w$, and a path $\mu$ from $u$ to $w$ would satisfy
$\mu=\mu w\in S$ by $(\#)$. Let $e$ be an edge on a path in $S$ and $v=s(e)$.
Then $v,r(e)\in S^{0}$ and there is a path $\mu$ from $r(e)$ to $v$. Thus
$d=e\mu$ is a closed path based at $v$, of length $\ell>0$, and by
Lemma~\ref{lem:ind}(c) the non-zero element $x=\mathbf{1}_{vS}$ satisfies
\[
x\cdot d=x .
\]
Let $\theta:N_{S}\rightarrow M$ be a homomorphism, where $M=\bigoplus
_{n\in\mathbb{Z}}M_{n}$ is graded, and put $a=\theta(x)$. Then $a\cdot d=a$.
Write $a=\sum_{n}a_{n}$ with $a_{n}\in M_{n}$. Since $d\in L_{\ell}$,
comparing components of degree $n+\ell$ gives $a_{n+\ell}=a_{n}\cdot d$ for
all $n$. If $a\neq0$, let $n_{0}$ be the least integer with $a_{n_{0}}\neq0$;
then $a_{n_{0}}=a_{n_{0}-\ell}\cdot d=0$, a contradiction. So $\theta(x)=0$.
By Lemma~\ref{lem:ind}(c),(d) and strong connectivity, $\mathbf{1}_{uS}%
=x\cdot\lambda$ for a path $\lambda$ from $v$ to $u\in S^{0}$, and
$\mathbf{1}_{\beta S}=\mathbf{1}_{r(\beta)S}\cdot\beta^{\ast}$; so
$N_{S}=x\cdot L_{K}(E)$ and $\theta=0$.
\end{proof}

\begin{corollary}
\label{cor:hom} Suppose $S$ is not of the form $\{w\}$, with $w$ a vertex. Then:

\begin{enumerate}
\item[(a)] $N_{S}$ is not a $\mathbb{Z}$-graded $L_{K}(E)$-module.

\item[(b)] $\operatorname{Hom}_{L_{K}(E)}(N_{S},L_{K}(E))=0$.

\item[(c)] $N_{S}$ is not a projective $L_{K}(E)$-module.
\end{enumerate}
\end{corollary}

\begin{proof}
(a) Otherwise the identity map of $N_{S}\neq0$ would contradict
Theorem~\ref{thm:hom}. (b) $L_{K}(E)$ is a graded right module over itself.
(c) Since $N_{S}$ is unitary, $N_{S}=\sum x\cdot v L_{K}(E)$, the sum taken
over all $x\in N_{S}$ and $v\in E^{0}$, so $N_{S}$ is a homomorphic image of a
direct sum $P$ of modules of the form $v L_{K}(E)$. Each $v L_{K}(E)$ is a
graded right ideal, so $P$ is a graded module. If $N_{S}$ were projective it
would be isomorphic to a direct summand of $P$, giving a non-zero homomorphism
$N_{S}\rightarrow P$.
\end{proof}

\begin{remark}
If $S=\{w\}$, then $N_{S}=K\delta_{w}$ is one-dimensional and graded.
\end{remark}

\subsection*{The annihilator}

\begin{proposition}
\label{prop:ann} Suppose that $S$ is an infinite set. Then the annihilator of
the simple module $N_{S}$ is $\operatorname{ann}(N_{S})=I(H,B_{H})$, a graded
primitive ideal.
\end{proposition}

\begin{proof}
Let $P=\operatorname{ann}(N_{S})$. By Lemma~\ref{lem:basic}, $H=E^{0}%
\backslash S^{0} $ is hereditary and saturated.

\emph{$I(H,B_{H})\subseteq P$.} Let $g\in N_{S}$, so $\operatorname{supp}%
(g)\subseteq S$. If $u\in H$, then $(g\cdot u)(y)$ is $g(y)$ if $s(y)=u$ and
$0$ otherwise; as no path in $S$ starts in $H$, $g\cdot u=0$. Let $u\in B_{H}%
$. For $y\in D$, by \eqref{eq:action}, $u^{H}\cdot y=y$ if $y=u$ or if $y$
begins with an edge $e\in s^{-1}(u)$ with $r(e)\in H$, and $u^{H}\cdot y=0$
otherwise. In the first case $y\notin S$, because $u\notin S$ by hypothesis;
in the second case $y\notin S$, because $r(e)\notin S^{0}$. Hence $(g\cdot
u^{H})(y)=g(u^{H}\cdot y)=0$ for all $y$. So $H\cup\{u^{H}:u\in B_{H}%
\}\subseteq P$.

\emph{$P\subseteq I(H,B_{H})$.} For $v\in S^{0}$ we have $\mathbf{1}_{vS}\cdot
v=\mathbf{1}_{vS}\neq0$, so $P\cap E^{0}=H$. Let $F$ be the graph of
Theorem~\ref{thm:known}(a), so $F^{0}=S^{0}$ and, by Lemma~\ref{lem:basic}(c),
$F^{1}$ is the set of edges of $E$ with source and range in $S^{0}$. Under the
isomorphism $L_{K}(E)/I(H,B_{H})\cong L_{K}(F)$, the ideal $P/I(H,B_{H})$
corresponds to an ideal $J$ of $L_{K}(F)$ containing no vertex. By
Theorem~\ref{thm:known}(b) it suffices to show that $F$ satisfies Condition (L).

Suppose $c$ is a cycle in $F$ without exits in $F$, and let $v\in c^{0}$. Let
$u\in S^{0}$. A path in $E$ from $v$ to $u$ has all its vertices in $S^{0}$
(Lemma~\ref{lem:basic}(c)), so it is a path in $F$; as $c$ has no exits in
$F$, this path runs along $c$ and $u\in c^{0}$. Hence $F^{0}=c^{0}$, and every
edge of $F$ is an edge of $c$. Every path in $S$ is a path in $F$. So the
infinite paths in $S$ are among the $|c|$ rotations of $c^{\infty}$. If $S$
contained a finite path, it would contain a singular vertex $w\in c^{0}$; as
$s_{S}(w)$ consists of the single edge of $c$ starting at $w$, $w$ would be an
infinite emitter in $B_{H}$, contrary to hypothesis. Thus $S$ is finite,
contrary to hypothesis. So $F$ satisfies Condition (L), $J=0$ and
$P=I(H,B_{H})$. This ideal is graded by Theorem~\ref{thm:known}(a), and it is
primitive because it is the annihilator of a simple module (see also
\cite{r-prime}).
\end{proof}

\begin{remark}
When $S$ is finite, $\operatorname{ann}(N_{S})$ need not be graded. Let
$S=\{c_{u}^{\infty}:u\in c^{0}\}$ be as in Lemma~\ref{lem:finiteS}(iv) and let
$v=s(c)$. As before, $I(H,B_{H})\subseteq\operatorname{ann}(N_{S})$. Now
$\mathbf{1}_{uS}\cdot(v-c)=0$ for all $u\in c^{0}$, so $v-c\in
\operatorname{ann}(N_{S})$. But $v-c\notin I(H,B_{H})$: this ideal is graded,
so it would contain the homogeneous component $v$ of $v-c$, whereas
$I(H,B_{H})\cap E^{0}=H$ and $v\notin H$.
\end{remark}

\subsection*{Finite presentation}

The next theorem describes when $N_{S}$ is finitely presented.

\begin{theorem}
\label{thm:fp} Let $N_{S}$ be the simple module described in
Theorem~\ref{thm:simple} (thus $S$ is an admissible set which contains no
breaking vertex of $H=E^{0}\setminus S^{0}$), where we assume that $S$ is not
the singleton set $\{w\}$ with $w$ a vertex. Then $N_{S}$ is a finitely
presented $L_{K}(E)$-module if and only if $S^{0}$ is finite and the number of
edges $e$ of $E$ with both $s(e)$ and $r(e)$ belonging to $S^{0}$ is finite.
\end{theorem}

\begin{proof}
As in the proof of Theorem~\ref{thm:hom}, $S$ contains a path of positive
length, some vertex $v\in S^{0}$ is the base of a closed path $d$, and
$N_{S}=x\cdot L_{K}(E)$ for $x=\mathbf{1}_{vS}$. Let $A$ denote the set of all
edges $e$ of $E$ with both $s(e)$ and $r(e)$ in $S^{0}$; thus $A$ is the union
of the sets $s_{S}(u)$, $u\in S^{0}$. By strong connectivity every vertex of
$S^{0}$ lies on a closed path, all of whose vertices are in $S^{0}$; so every
vertex of $S^{0}$ is the source of at least one edge in $A$. In particular,
$S^{0}$ is finite if $A$ is finite.

\textit{Necessity.} By Lemma~\ref{lem:source}, the tail-equivalence class
$[d^{\infty}]$ is an admissible set of infinite paths with $[d^{\infty}%
]^{0}=S^{0}$, so $N_{S}\cong N_{[d^{\infty}]}$ by Proposition~\ref{prop:indep}%
. Hence, in proving necessity, we may assume that $S=[d^{\infty}]$; in
particular $d^{\infty}\in vS$.

Suppose that $N_{S}$ is finitely presented and assume, by way of
contradiction, that the set $A$ is infinite. The map $vL_{K}(E)\rightarrow
N_{S}$, $a\mapsto x\cdot a$, is surjective and $vL_{K}(E)$ is a finitely
generated projective module. So its kernel $I=\{a\in vL_{K}(E):x\cdot a=0\}$
is finitely generated, say $I=\sum_{j=1}^{n}a_{j}L_{K}(E)$. Write each $a_{j}
$ as a $K$-linear combination of monomials $\alpha\beta^{\ast}$, and choose an
edge $e_{0}\in A$ which occurs neither in $d$ nor in any of the paths
$\alpha,\beta$ appearing in these finitely many monomials. Let
\[
Y=\{p\in vS:e_{0}\text{ does not occur in }p\},\qquad y=\mathbf{1}_{Y}\in
K^{D}.
\]
Then $y\neq0$ since $d^{\infty}\in Y$.

We show that $y\cdot a_{j}=0$ for all $j$. Let $\alpha\beta^{\ast}$ be a
monomial in which $e_{0}$ does not occur, and let $p\in D$. If $\beta
\not \preceq p$, then $(x\cdot\alpha\beta^{\ast})(p)=0=(y\cdot\alpha
\beta^{\ast})(p)$ by \eqref{eq:right}. If $p=\beta p^{\prime}$, then
$(x\cdot\alpha\beta^{\ast})(p)=\mathbf{1}_{vS}(\alpha p^{\prime})$ and
$(y\cdot\alpha\beta^{\ast})(p)=\mathbf{1}_{Y}(\alpha p^{\prime})$.

Suppose first that $e_{0}$ does not occur in $p$. If $p=\beta p^{\prime}$,
then $e_{0}$ does not occur in $\alpha p^{\prime}$, so $\alpha p^{\prime}$
belongs to $Y$ if and only if it belongs to $vS$, and $\mathbf{1}_{Y}(\alpha
p^{\prime})=\mathbf{1}_{vS}(\alpha p^{\prime})$. Thus $(y\cdot\alpha
\beta^{\ast})(p)=(x\cdot\alpha\beta^{\ast})(p)$ for every monomial
$\alpha\beta^{\ast}$ of $a_{j}$, and by linearity $(y\cdot a_{j})(p)=(x\cdot
a_{j})(p)=0$, since $a_{j}\in I$.

Suppose next that $e_{0}$ occurs in $p$. If $p=\beta p^{\prime}$, then $e_{0}$
occurs in $p^{\prime}$, because it does not occur in $\beta$; so $\alpha
p^{\prime}\notin Y$ and $(y\cdot\alpha\beta^{\ast})(p)=0$. Thus $(y\cdot
\alpha\beta^{\ast})(p)=0$ for every monomial of $a_{j}$, and $(y\cdot
a_{j})(p)=0$. Hence $y\cdot a_{j}=0$ for all $j$, and therefore $y\cdot I=0$.

On the other hand, since $s(e_{0}),r(e_{0})\in S^{0}$, there are paths
$\lambda$ from $v$ to $s(e_{0})$ and $\mu$ from $r(e_{0})$ to $v$. Put
$b=\lambda e_{0}\mu-v\in vL_{K}(E)$. By Lemma~\ref{lem:ind}(c), $x\cdot\lambda
e_{0}\mu=\mathbf{1}_{vS}=x\cdot v$, so $x\cdot b=0$ and $b\in I$. But, for
every $p\in D$, $(y\cdot\lambda e_{0}\mu)(p)$ is either $0$ or $\mathbf{1}%
_{Y}(\lambda e_{0}\mu p)$, and the latter is $0$ because $e_{0}$ occurs in
$\lambda e_{0}\mu p$; while $y\cdot v=y$. So $y\cdot b=-y\neq0$, contradicting
$y\cdot I=0$. Therefore $A$ is a finite set, and so is $S^{0}$.

\textit{Sufficiency.} Suppose that the set $A$ is finite; then $S^{0}$ is
finite. We first show that $S\subseteq D_{\infty}$. Otherwise $S$ contains a
finite path and hence a vertex $w$ (Lemma~\ref{lem:basic}(a)), which is
singular. Since $S\neq\{w\}$, Lemma~\ref{lem:single} gives $s_{S}%
(w)\neq\emptyset$; so $w$ is not a sink, hence $w$ is an infinite emitter. As
$s_{S}(w)=s^{-1}(w)\cap r^{-1}(E^{0}\setminus H)$ is contained in $A$, it is
finite, and so $w\in B_{H}$, contrary to the hypothesis.

Let $\beta$ be a finite path with $r(\beta)=v\in S^{0}$. As $S\subseteq
D_{\infty}$, every $p\in\beta S$ is an infinite path, so $p=\beta ep^{\prime}$
for a unique edge $e$ with $s(e)=v$, and $r(e)\in S^{0}$ since $r(e)$ lies on
$p$; thus $e\in s_{S}(v)$. Hence $\beta S$ is the disjoint union of the sets
$\beta eS$, $e\in s_{S}(v)$, and
\begin{equation}
\mathbf{1}_{\beta S}=\sum_{e\in s_{S}(v)}\mathbf{1}_{\beta eS}\qquad
(v=r(\beta)\in S^{0}). \label{eq:split}%
\end{equation}

Since $S^{0}$ is finite, $\varepsilon=\sum_{v\in S^{0}}v$ is an idempotent of
$L_{K}(E)$, and $P=\varepsilon L_{K}(E)=\bigoplus_{v\in S^{0}}vL_{K}(E)$ is a
finitely generated projective right $L_{K}(E)$-module. Let $R$ be the
submodule of $P$ generated by the finitely many elements
\[
e-r(e)\ \ (e\in A),\qquad v-\sum_{e\in s_{S}(v)}ee^{\ast}\ \ (v\in S^{0}),
\]
and put $M=P/R$, a finitely presented module. Let $x_{v}$ be the image of $v$
in $M$. We wish to show that $M\cong N_{S}$.

Define $\pi:P\rightarrow N_{S}$ by $\pi(a)=\mathbf{1}_{S}\cdot a$. Since
$\mathbf{1}_{S}\cdot u=\mathbf{1}_{uS}=0$ for $u\notin S^{0}$
(Lemma~\ref{lem:ind}(b)), we have $\mathbf{1}_{S}\cdot a=\mathbf{1}_{S}%
\cdot\varepsilon a$ for every $a\in L_{K}(E)$, so $\pi$ is surjective. We next
show that $\pi(R)=0$.

By Lemma~\ref{lem:ind}(b),(c), $\mathbf{1}_{S}\cdot e=\mathbf{1}_{s(e)S}\cdot
e=\mathbf{1}_{r(e)S}=\mathbf{1}_{S}\cdot r(e)$ for $e\in A$, and so
$\pi(e-r(e))=\mathbf{1}_{S}\cdot(e-r(e))=0$.

Let $v\in S^{0}$. (If $v$ is a regular vertex, then by (CK-2) the element
$v-\sum_{e\in s_{S}(v)}ee^{\ast}$ equals $\sum_{e\in s^{-1}(v),\,r(e)\in
H}ee^{\ast}$; if $v$ is an infinite emitter, then $s_{S}(v)$ is finite and
non-empty, so $v\in B_{H}$ and this element is $v^{H}$.) By
Lemma~\ref{lem:ind} and \eqref{eq:split} (with $\beta=v$),
\[
\mathbf{1}_{S}\cdot\Bigl(v-\sum_{e\in s_{S}(v)}ee^{\ast}\Bigr)=\mathbf{1}%
_{vS}-\sum_{e\in s_{S}(v)}\mathbf{1}_{eS}=0.
\]
So $\pi$ vanishes on $R$ and induces a surjection $\bar{\pi}:M\rightarrow
N_{S}$ with $\bar{\pi}(x_{v})=\mathbf{1}_{vS}$. For a finite path $\beta$ with
$r(\beta)\in S^{0}$ put $m_{\beta}=x_{r(\beta)}\beta^{\ast}\in M$; then
$\bar{\pi}(m_{\beta})=\mathbf{1}_{\beta S}$ by Lemma~\ref{lem:ind}(d).

\emph{We claim that for every $v\in S^{0}$,}
\begin{equation}
x_{v}=\sum_{e\in s_{S}(v)}x_{r(e)}e^{\ast}\qquad\text{and}\qquad
x_{v}f=0\ \text{ for every edge }f\notin s_{S}(v). \label{eq:xv}%
\end{equation}
Indeed, the defining relations of $M$ give $x_{v}=\sum_{e\in s_{S}(v)}%
x_{v}ee^{\ast}$ and $x_{v}e=x_{r(e)}$ for $e\in A$, which is the first
statement. Let $f\notin s_{S}(v)$ be an edge. If $s(f)\neq v$, then
$x_{v}f=x_{v}vf=0$. If $s(f)=v$, then $x_{v}f=\sum_{e\in s_{S}(v)}%
x_{v}e(e^{\ast}f)=0$ by (CK-1).

\emph{We claim that $\bar{\pi}$ is injective.} By \eqref{eq:xv} and the
relations $x_{s(e)}e=x_{r(e)}$ ($e\in A$), we have $x_{v}\alpha=x_{r(\alpha)}$
if $\alpha$ is a path starting at $v$ all of whose vertices lie in $S^{0}$,
and $x_{v}\alpha=0$ for every other finite path $\alpha$. Hence $x_{v}%
\alpha\beta^{\ast}$ is either $0$ or $m_{\beta}$, and $M$ is spanned by the
elements $m_{\beta}$. Again by \eqref{eq:xv}, $m_{\beta}=\sum_{e\in
s_{S}(r(\beta))}m_{\beta e}$, and by induction $m_{\beta}$ is the sum of the
elements $m_{\beta\mu}$, where $\mu$ runs over the paths of any given length
starting at $r(\beta)$ all of whose vertices lie in $S^{0}$. Now let $z\in
\ker\bar{\pi}$. By the above we can write $z=\sum_{j}c_{j}m_{\beta_{j}}$ with
$c_{j}\in K$ and the $\beta_{j}$ distinct paths of one and the same length
with $r(\beta_{j})\in S^{0}$. Then $0=\bar{\pi}(z)=\sum_{j}c_{j}%
\mathbf{1}_{\beta_{j}S}$. The sets $\beta_{j}S$ are pairwise disjoint and
non-empty (Lemma~\ref{lem:basic}(e)), so all $c_{j}=0$ and $z=0$.

Thus $\bar\pi:M\rightarrow N_{S}$ is an isomorphism, and $N_{S}\cong P/R$ is
finitely presented.
\end{proof}

\begin{example}
For the graph $R_{n}$ ($n\geq2$) of Example~\ref{ex:rose} and any admissible
set $S$, Theorem~\ref{thm:fp} gives
\[
N_{S}\cong L_{K}(R_{n})\Big/\sum_{i=1}^{n}(e_{i}-v)L_{K}(R_{n}).
\]
On the other hand, the simple module of Example~\ref{ex:emitter} is not
finitely presented, since infinitely many loops are based at $w$.
\end{example}

\begin{remark}
\label{rem:fpbreak} In Theorem~\ref{thm:fp} one must distinguish between $S$,
a set of paths, and $S^{0}$, the set of vertices lying on those paths. A
breaking vertex $u$ of $H$ always belongs to $S^{0}$; the hypothesis $S\cap
B_{H}=\emptyset$ only says that the path $u$ of length $0$ is not an element
of $S$. For example, let $E$ have a vertex $w$, two loops $e,f$ based at $w$,
and infinitely many edges from $w$ to sinks, and let $S$ be the set of all
infinite paths in $E$. Then $S^{0}=\{w\}$, $s_{S}(w)=\{e,f\}$ and $w\in B_{H}%
$, while $w\notin S$. By Theorem~\ref{thm:fp}, $N_{S}$ is finitely presented,
although $w$ is an infinite emitter. Under the hypotheses of the theorem,
every vertex $u$ which is an element of $S$ is an infinite emitter with
$s_{S}(u)$ infinite (it is not a sink by Lemma~\ref{lem:single}, and it is not
in $B_{H}$); in particular, no such vertex exists when the finiteness
condition of the theorem holds.
\end{remark}

It remains to consider the case $S=\{w\}$.

\begin{proposition}
\label{prop:acyclic} Let $S=\{w\}$ be admissible. If $w$ is an isolated
vertex, then $N_{S}\cong wL_{K}(E)$ is projective and finitely presented. If
$w$ is an infinite emitter, then $N_{S}$ is not finitely presented.
\end{proposition}

\begin{proof}
Here $w$ is a source and $N_{S}=K\delta_{w}$. Let $\pi:wL_{K}(E)\rightarrow
N_{S}$, $\pi(a)=\delta_{w}\cdot a$. By \eqref{eq:right}, for finite paths
$\alpha,\beta$ with $s(\alpha)=w$ and $r(\alpha)=r(\beta)$, $\delta_{w}%
\cdot\alpha\beta^{\ast}$ is $\delta_{\beta}$ if $\alpha=w$ and $0$ otherwise;
and $\alpha=w$ forces $r(\beta)=w$, hence $\beta=w$ because $w$ is a source.
So $\pi(w)=\delta_{w}$ and $\pi(\alpha\beta^{\ast})=0$ for all other such
monomials, each of which lies in $eL_{K}(E)$ for the first edge $e$ of
$\alpha$. Therefore $wL_{K}(E)=Kw\oplus\ker\pi$ and
\[
\ker\pi=\sum_{e\in s^{-1}(w)}eL_{K}(E)=\bigoplus_{e\in s^{-1}(w)}ee^{\ast
}L_{K}(E),
\]
the sum being direct because the idempotents $ee^{\ast}$ are pairwise
orthogonal. If $w$ is isolated, $\ker\pi=0$ and $N_{S}\cong wL_{K}(E)$. If $w$
is an infinite emitter, $\ker\pi$ is a direct sum of infinitely many non-zero
right ideals, so it is not finitely generated, and $N_{S}$ is not finitely presented.
\end{proof}

\begin{corollary}
\label{cor:allfp}

\begin{enumerate}
\item[(a)] If $E$ is a finite graph, then every simple module $N_{S}$ is
finitely presented.

\item[(b)] If $E$ is a row-finite graph, then every simple module $N_{S}$ is
finitely presented if and only if $S^{0}$ is finite for every admissible set
$S$ which is not a single vertex.
\end{enumerate}
\end{corollary}

\begin{proof}
Let $E$ be row-finite, so that $E$ has no infinite emitters. If $S=\{w\}$ is
admissible, then $w$ is a singular vertex which is a source, hence an isolated
vertex, and $N_{S}$ is finitely presented by Proposition~\ref{prop:acyclic}.
Now let $S$ be an admissible set which is not a single vertex. Since $E$ has
no breaking vertices, $N_{S}$ is simple (Corollary~\ref{Finite case}). Let $A$
be the set of edges with both source and range in $S^{0}$. Every vertex of
$S^{0}$ is the source of at least one edge in $A$ (Lemma~\ref{lem:source}(b)),
and each vertex emits only finitely many edges, so $A$ is finite if and only
if $S^{0}$ is finite. By Theorem~\ref{thm:fp}, $N_{S}$ is finitely presented
if and only if $S^{0}$ is finite. This proves (b).

(a) A finite graph is row-finite and $S^{0}\subseteq E^{0}$ is always finite,
so (a) follows from (b).
\end{proof}

\begin{remark}
By contrast, when $E$ is a finite graph, every simple $L_{K}(E)$-module is
finitely presented if and only if every vertex of $E$ is the base of at most
one cycle \cite[Theorem 1.1]{ar}. So, over a finite graph, the simple modules
$N_{S}$ are always finitely presented, while other simple modules need not be.
\end{remark}

\section{A left action of $L_{K}(E)$ on $K^{D}$}

\label{sec:left}

Throughout this section $S$ is an admissible set of paths, $H=E^{0}\setminus
S^{0}$, and $N_{S}=\mathbf{1}_{S}\cdot L_{K}(E)$ as before. We now let
$L_{K}(E)$ act on $K^{D}$ \emph{on the left}, and ask when the left module
generated by $\mathbf{1}_{S}$ is simple. The motivation for this is our desire
to compare the simple right module $N_{S}$ with the Chen simple modules, all
of which are defined as left $L_{K}(E)$-modules. When $S$ is infinite,
$\mathbf{1}_{S}$ is an element of the dual space $K^{D}$ and not of $KD$
itself, so a left action of $L_{K}(E)$ on $\mathbf{1}_{S}$ has to be built
from the right action $\cdot$ already defined on $K^{D}$; twisting it through
the involution $\ast$ is a natural way to do this without leaving $K^{D}$. We
show that the resulting left module $L_{K}(E)\mathbf{1}_{S}$ is simple in
exactly the same cases as $N_{S}$, and we compare it with the Chen simple left modules.

For $a\in L_{K}(E)$ and $f\in K^{D}$ put
\[
a\circ f=f\cdot a^{\ast}.
\]

\begin{lemma}
\label{lem:circ} The rule above makes $K^{D}$ a left $L_{K}(E)$-module.
Moreover a subspace $W\subseteq K^{D}$ satisfies $L_{K}(E)\circ W\subseteq W$
if and only if $W\cdot L_{K}(E)\subseteq W$; thus $W$ is a submodule of
$(K^{D},\circ)$ if and only if it is a submodule of $(K^{D},\cdot)$.
\end{lemma}

\begin{proof}
For $a,b\in L_{K}(E)$ and $f\in K^{D}$, using $(ab)^{\ast}=b^{\ast}a^{\ast}$
and associativity of the right action $\cdot$,
\[
a\circ(b\circ f)=a\circ(f\cdot b^{\ast})=(f\cdot b^{\ast})\cdot a^{\ast
}=f\cdot(b^{\ast}a^{\ast})=f\cdot(ab)^{\ast}=(ab)\circ f,
\]
so $\circ$ is a left action. Since $a\mapsto a^{\ast}$ is a bijection of
$L_{K}(E)$ onto itself, the statement ``$a\circ f\in W$ for every $a\in
L_{K}(E)$ and $f\in W$'' is equivalent to ``$f\cdot b\in W$ for every $b\in
L_{K}(E)$ and $f\in W$''.
\end{proof}

\begin{lemma}
\label{lem:circ1S} For finite paths $\alpha,\beta$ with $r(\alpha)=r(\beta)$,
$(\alpha\beta^{\ast})\circ\mathbf{1}_{S}=\mathbf{1}_{\alpha S}$. Consequently
$L_{K}(E)\circ\mathbf{1}_{S}=N_{S}$, as subsets of $K^{D}$.
\end{lemma}

\begin{proof}
Since $(\alpha\beta^{\ast})^{\ast}=\beta\alpha^{\ast}$,
\[
(\alpha\beta^{\ast})\circ\mathbf{1}_{S}=\mathbf{1}_{S}\cdot(\alpha\beta^{\ast
})^{\ast}=\mathbf{1}_{S}\cdot\beta\alpha^{\ast}=\mathbf{1}_{\alpha S}%
\]
by Lemma~\ref{lem:ind}, applied with the roles of its two paths exchanged.
Taking $K$-linear combinations, $L_{K}(E)\circ\mathbf{1}_{S}%
=\operatorname{span}_{K}\{\mathbf{1}_{\alpha S}:\alpha\text{ a finite path
with }r(\alpha)\in S^{0}\}$, which is $N_{S}$ by Lemma~\ref{lem:ind}(a).
\end{proof}

\begin{proposition}
\label{prop:leftstar} The left module $(N_{S},\circ)$ is simple if and only if
$S$ contains no breaking vertex of $H$, that is, if and only if the right
module $N_{S}=(N_{S},\cdot)$ is simple.
\end{proposition}

\begin{proof}
Immediate from Lemma~\ref{lem:circ1S}, Lemma~\ref{lem:circ} and
Theorem~\ref{thm:simple}: $(N_{S},\circ)$ and $(N_{S},\cdot)$ have the same
subspaces as submodules.
\end{proof}

\begin{remark}
Thus the left module $L_{K}(E)\mathbf{1}_{S}=L_{K}(E)\circ\mathbf{1}_{S}$ is
simple in exactly the same cases as the right module $\mathbf{1}_{S}\cdot
L_{K}(E)$, with the same necessary and sufficient condition: $S$ contains no
breaking vertex of $H$. In particular Example~\ref{ex:toeplitz}, read through
$\circ$, already exhibits a non-simple instance of this left module.
\end{remark}

\subsection*{Comparison with the Chen simple modules}

We now compare the simple left module $L_{K}(E)\mathbf{1}_{S}=(N_{S},\circ)$
of Proposition~\ref{prop:leftstar} with the Chen simple left modules. We
recall these from \cite{chen}, \cite{ar}, \cite{r-simple}:

\begin{enumerate}
\item[(C1)] $V_{[p]}$, for an infinite path $p$: the $K$-space with basis the
tail-equivalence class $[p]$, with $v\cdot q$, $e\cdot q$ and $e^{\ast}\cdot
q$ defined as for $KD$ above.

\item[(C2)] The modules attached to a sink or to an infinite emitter $w$:
$K$-spaces having a basis $B$ of finite paths ending at $w$, on which, for
each edge $e$ and each $q\in B$, $e^{\ast}\cdot q$ is either $0$ or the path
$q^{\prime}$ with $q=eq^{\prime}$.

\item[(C3)] $V^{f}_{c}=L_{K}(E) v/L_{K}(E) f(c)$, where $c$ is a cycle based
at $v$ and $f(t)\in K[t]$ is a polynomial with $f(0)\neq0$.
\end{enumerate}

\begin{proposition}
\label{prop:poly} Let $S$ be an infinite admissible set containing no breaking
vertex of $H$. Let $c$ be a closed path in $E$ based at a vertex $v$, let
$f(t)=\sum_{j=0}^{n}a_{j}t^{j}\in K[t]$ with $a_{0}\neq0$, and put
$f(c)=a_{0}v+\sum_{j=1}^{n}a_{j}c^{j}\in L_{K}(E)$. If $x\in N_{S}$ satisfies
$v\circ x=x$ and $f(c)\circ x=0$, then $x=0$. In particular, no non-zero $x\in
N_{S}$ satisfies $c\circ x=x$.
\end{proposition}

\begin{proof}
Since $f(c)^{\ast}=a_{0}v+\sum_{j=1}^{n}a_{j}(c^{\ast})^{j}=:f(c^{\ast})$, the
definition of $\circ$ shows that the hypotheses say $x\cdot v=x$ and $x\cdot
f(c^{\ast})=0$. Since $x=x\cdot v$, $x(p)=0$ whenever $s(p)\neq v$. For $p\in
D$ with $s(p)=v $ let $k(p)$ be the supremum of the integers $k\geq0$ with
$c^{k}\preceq p$. Then $k(p)=\infty$ only if $p=c^{\infty}$. We show by
induction on $k(p)$ that $x(p)=0$ whenever $k(p)<\infty$. By \eqref{eq:right},
$(x\cdot(c^{\ast})^{j})(p)=x(p_{j})$ if $p=c^{j}p_{j}$, and $0$ if
$c^{j}\not \preceq p$. Hence
\[
0=(x\cdot f(c^{\ast}))(p)=a_{0}x(p)+\sum_{j=1}^{\min(n,k(p))}a_{j}x(p_{j}).
\]
For each $j$ in this sum we have $s(p_{j})=v$ and $k(p_{j})=k(p)-j<k(p)$, so
$x(p_{j})=0$ by induction, and then $x(p)=0$ because $a_{0}\neq0$. Therefore
$\operatorname{supp}(x)\subseteq\{c^{\infty}\}$, and $x=0$ by
Corollary~\ref{cor:support}.

For the last statement take $f(t)=1-t$, so that $f(c)=v-c$. If $c\circ x=x$,
then $v\circ x=v\circ(c\circ x)=(vc)\circ x=c\circ x=x$, and $f(c)\circ
x=x-x=0$; so $x=0$.
\end{proof}

\begin{theorem}
\label{thm:chen} Let $S$ be an infinite admissible set containing no breaking
vertex of $H$. Then the simple left module $L_{K}(E)\mathbf{1}_{S}%
=(N_{S},\circ)$ is not isomorphic to any Chen simple left module $V$ of type
\textup{(C1)}, \textup{(C2)} or \textup{(C3)}.
\end{theorem}

\begin{proof}
By the proof of Theorem~\ref{thm:hom} there are $v\in S^{0}$, a closed path
$d$ based at $v$ of length $\ell>0$, and a non-zero $x\in N_{S}$ with $x\cdot
v=x$ and $x\cdot d=x$, that is, $v\circ x=x$ and $d^{\ast}\circ x=x$. Suppose
$\varphi:(N_{S},\circ)\rightarrow V$ is an isomorphism of left modules and let
$y=\varphi(x)\neq0$. Then $vy=y$ and $d^{\ast}y=y$ in $V$.

(C1) Write $y=\sum_{j=1}^{m}k_{j}q_{j}$ with the $q_{j}\in[p]$ distinct and
all $k_{j}\neq0$. Then $d^{\ast}y=\sum\{k_{j}q_{j}^{\prime}:q_{j}%
=dq_{j}^{\prime}\}$, and the paths $q_{j}^{\prime}$ are distinct. From
$d^{\ast}y=y$ it follows that every $q_{j}$ begins with $d$ and that $q\mapsto
q^{\prime}$ is a bijection of the finite set $Q=\{q_{1},\dots,q_{m}\}$ onto
itself. Hence each $q\in Q$ begins with $d^{k}$ for every $k$, so
$Q=\{d^{\infty}\}$ and $y=k_{1}d^{\infty}$. Then $dy=y$ in $V$, so $d\circ
x=\varphi^{-1}(dy)=x$. This contradicts Proposition~\ref{prop:poly}.

(C2) Write $y=\sum_{j}k_{j}q_{j}$ with $q_{j}\in B$ distinct and $k_{j}\neq0$.
Each $d^{\ast}q_{j}$ is $0$ or a path of length $|q_{j}|-\ell$. So the maximal
length of a path occurring in $d^{\ast}y$ is smaller than the maximal length
of a path occurring in $y$, and $d^{\ast}y=y$ is impossible.

(C3) Let $y_{0}$ be the image of $v$ in $V=V^{f}_{c}$. Then $vy_{0}=y_{0}$ and
$f(c)y_{0}=0$. The non-zero element $x_{0}=\varphi^{-1}(y_{0})\in N_{S}$
therefore satisfies $v\circ x_{0}=x_{0}$ and $f(c)\circ x_{0}=0$,
contradicting Proposition~\ref{prop:poly}.
\end{proof}

\begin{remark}
(1) The modules of type (C1) with $p$ irrational and those of type (C2) are
graded \cite{hr}. If $V=\bigoplus_{n}V_{n}$ is such a module, then $V$ is a
graded right module under $y\cdot a=a^{\ast}y$, with $V_{n}$ placed in degree
$-n$, and an isomorphism $(N_{S},\circ)\cong V$ would be an isomorphism of
right modules from $N_{S}$ onto it. So for these modules the conclusion of
Theorem~\ref{thm:chen} also follows from Corollary~\ref{cor:hom}(a).

(2) The hypothesis that $S$ is infinite is needed: by Example~\ref{ex:loop},
for the finite set $S=\{c_{u}^{\infty}:u\in c^{0}\}$ the element
$x=\mathbf{1}_{s(c)S}$ satisfies $c\circ x=x$.

(3) Theorem~\ref{thm:chen} is consistent with \cite[Theorem 1.1]{ar},
according to which every simple module over $L_{K}(E)$ is a Chen module when
$E $ is a finite graph in which every vertex is the base of at most one cycle.
In such a graph every admissible set $S$ is finite. Indeed, if $S$ is not a
single vertex, then $S\subseteq D_{\infty}$ (Lemma~\ref{lem:single}) and the
subgraph $F$ on $S^{0}$ contains a cycle $c$. If $F$ had an edge not on $c$,
then by strong connectivity it would have such an edge $f$ with $s(f)\in
c^{0}$; a shortest path $\mu$ in $F$ from $r(f)$ to $s(f)$ would give a cycle
$f\mu$ based at $s(f)$ and different from the rotation of $c$ based at $s(f)$.
So $F$ consists of $c$ alone, and $S$, whose elements are infinite paths in
$F$, is finite.
\end{remark}

\subsection*{The simple module $N_{S}$ when $E$ is the rose graph $R_{n}$}

T.G. Nam asked whether the simple module $N_{S}$ is isomorphic to any of the
simple modules over the Leavitt path algebra $L_{K}(R_{n})$ constructed in his
joint paper \cite{nsv}. Here we provide a complete answer. Throughout this
subsection $n\geq2$, $E=R_{n}$ is the rose graph with one vertex $v$ and $n$
loops $e_{1},\dots,e_{n}$, and $L=L_{K}(R_{n})$. Then $v$ is the identity of
$L$, $v$ is a regular vertex and $D=D_{\infty}$. By Example~\ref{ex:rose}, the
admissible sets are exactly the non-empty unions of tail-equivalence classes
of infinite paths, and each of them satisfies $S^{0}=\{v\}$. Every admissible
set $S$ is infinite by Lemma~\ref{lem:finiteS}, because each loop $e_{i}$ is
entered by the other loops, and $R_{n}$ has no breaking vertices. Hence all
the modules $N_{S}$ are simple (Corollary~\ref{cor:infinite}) and isomorphic
to one another (Proposition~\ref{prop:indep}); they are finitely presented
(Theorem~\ref{thm:fp}) and neither graded nor projective
(Corollary~\ref{cor:hom}). Since $S^{0}$ is finite, $\mathbf{1}_{S}%
=\mathbf{1}_{vS}$ belongs to $N_{S}$ and $N_{S}=\mathbf{1}_{S}\cdot L$. We
first show that $N_{S}$ and $L\mathbf{1}_{S}=(N_{S},\circ)$ are presented by a
single relation. Put
\[
\varepsilon=\sum_{k=1}^{n}e_{k}\in L,\qquad\text{so that}\qquad\varepsilon
^{\ast}=\sum_{k=1}^{n}e_{k}^{\ast}.
\]

\begin{proposition}
\label{prop:rose} Let $S$ be an admissible set of paths in $R_{n}$.

\begin{enumerate}
\item[(a)] $\sum_{i=1}^{n}(e_{i}-v)L=(v-\varepsilon^{\ast})L$, and
$a+(v-\varepsilon^{\ast})L\mapsto\mathbf{1}_{S}\cdot a$ is an isomorphism of
right $L$-modules from $L/(v-\varepsilon^{\ast})L$ onto $N_{S}$.

\item[(b)] $\sum_{i=1}^{n}L(e_{i}^{\ast}-v)=L(v-\varepsilon)$, and
$a+L(v-\varepsilon)\mapsto a\circ\mathbf{1}_{S}$ is an isomorphism of left
$L$-modules from $L/L(v-\varepsilon)$ onto $L\mathbf{1}_{S}=(N_{S},\circ)$.
\end{enumerate}
\end{proposition}

\begin{proof}
(a) Here $S^{0}=\{v\}$ and the set $A$ of edges with source and range in
$S^{0}$ is $\{e_{1},\dots,e_{n}\}=s_{S}(v)$. So in the proof of
Theorem~\ref{thm:fp} we have $P=vL=L$, the generators $v-\sum_{i}e_{i}%
e_{i}^{\ast}$ of $R$ are $0$ by (CK-2), and the map $\bar{\pi}$ is an
isomorphism $L/\sum_{i}(e_{i}-v)L\rightarrow N_{S}$, $a+\sum_{i}%
(e_{i}-v)L\mapsto\mathbf{1}_{S}\cdot a$. By (CK-2) and (CK-1),
\[
\sum_{k=1}^{n}(e_{k}-v)e_{k}^{\ast}=v-\varepsilon^{\ast}\qquad\text{and}%
\qquad(v-\varepsilon^{\ast})e_{i}=e_{i}-v\quad(1\leq i\leq n),
\]
so the right ideals $\sum_{i}(e_{i}-v)L$ and $(v-\varepsilon^{\ast})L$ coincide.

(b) Applying the involution to (a), and using $(ab)^{\ast}=b^{\ast}a^{\ast}$,
we get
\[
\sum_{i=1}^{n}L(e_{i}^{\ast}-v)=\Bigl(\sum_{i=1}^{n}(e_{i}-v)L\Bigr)^{\ast
}=\bigl((v-\varepsilon^{\ast})L\bigr)^{\ast}=L(v-\varepsilon).
\]
The map $L\rightarrow(N_{S},\circ)$, $a\mapsto a\circ\mathbf{1}_{S}%
=\mathbf{1}_{S}\cdot a^{\ast}$, is a homomorphism of left modules, it is
surjective by Lemma~\ref{lem:circ1S}, and by (a) its kernel is $\{a\in
L:a^{\ast}\in(v-\varepsilon^{\ast})L\}=L(v-\varepsilon)$.
\end{proof}

Nam, Srivastava and Vien \cite{nsv} constructed several classes of simple left
$L$-modules, by twisting Chen modules through automorphisms of $L$, and showed
that the modules in their list \cite[Theorem 4.7]{nsv} are pairwise
non-isomorphic. We now describe which of these simple modules is isomorphic to
$L\mathbf{1}_{S}=(N_{S},\circ)$. We first recall the notation of \cite{nsv}.

For $P=(p_{ij})\in GL_{n}(K)$ with $P^{-1}=(p_{ij}^{\prime})$ there is a
graded automorphism $\varphi_{P}$ of $L$ with
\[
\varphi_{P}(v)=v,\qquad\varphi_{P}(e_{i})=\sum_{k=1}^{n}e_{k}p_{ki}%
,\qquad\varphi_{P}(e_{i}^{\ast})=\sum_{k=1}^{n}p_{ik}^{\prime}e_{k}^{\ast
}\qquad(1\leq i\leq n)
\]
\cite[Proposition 2.7 and Corollary 2.8]{nsv}. If $M$ is a left $L$-module and
$\psi$ is an automorphism of $L$, then $M^{\psi}$ denotes the left $L$-module
with the same underlying space and the action $a\cdot m=\psi(a)m$. Clearly
$M^{\psi}$ is simple if and only if $M$ is simple, and $(M^{\psi})^{\psi^{-1}%
}=M$. For an infinite path $\alpha$ put $V_{[\alpha]}^{P}=(V_{[\alpha
]})^{\varphi_{P}^{-1}}$. A closed path $c$ is \emph{simple} if $c\neq d^{k}$
for every closed path $d$ and every $k\geq2$, and $\mathrm{SCP}(R_{n})$
denotes the set of simple closed paths in $R_{n}$. Let $\mathrm{Irr}(K[x])$ be
the set of irreducible polynomials of the form $f=1-a_{1}x-\cdots-a_{m}x^{m}$.
For $c\in\mathrm{SCP}(R_{n})$ and $f\in\mathrm{Irr}(K[x])$, $S_{\Pi_{c}}^{f}$
denotes the simple module of \cite{an}, generated by an element $z$ subject to
$f(c)z=0$, that is, $S_{\Pi_{c}}^{f}\cong L/Lf(c)$; the modules $S_{d}^{f,p}$
are the twists of the modules $S_{d}^{f}$ by the Anick type automorphisms of
\cite{kn}; and $(R_{n})_{\mathrm{irr-eeri}}^{\infty}$ is the set of irrational
infinite paths in which every loop occurs infinitely many times. The list in
\cite[Theorem 4.7]{nsv} consists of the following simple left $L$-modules:

\begin{enumerate}
\item[(1)] $V_{[\alpha]}$, where $\alpha$ is an irrational infinite path;

\item[(2)] $S_{\Pi_{c}}^{f}$, where $c\in\mathrm{SCP}(R_{n})$ and
$f\in\mathrm{Irr}(K[x])$;

\item[(3)] $S_{d}^{f,p}$, where $f\in\mathrm{Irr}(K[x])$ has degree at least
$2$ and $d$, $p$ are as in \cite[Theorem 4.7(3)]{nsv};

\item[(4)] $V_{[\alpha]}^{P}$, where $\alpha\in(R_{n})_{\mathrm{irr-eeri}%
}^{\infty}$ and $I_{n}\neq P\in GL_{n}(K)$;

\item[(5)] $V_{[c^{\infty}]}^{P}$, where $c\in\mathrm{SCP}(R_{n})$ and
$I_{n}\neq P\in GL_{n}(K)$.
\end{enumerate}

We shall use the following facts from \cite{nsv} and \cite{kn}: for
$c\in\mathrm{SCP}(R_{n})$ and $P\in GL_{n}(K)$, the module $V_{[c^{\infty}%
]}^{P}$ is simple, $\operatorname{End}_{L}(V_{[c^{\infty}]}^{P})\cong K$ and
\begin{equation}
V_{[c^{\infty}]}^{P}\cong L/L(v-\varphi_{P}(c)) \label{eq:nsvc}%
\end{equation}
\cite[Theorem 4.5]{nsv}; the modules in (4) are not finitely presented
\cite[Theorem 4.2(3)]{nsv}; and $\operatorname{End}_{L}(S_{d}^{f,p})\cong
K[x]/K[x]f(x)$ \cite[Theorem 3.6(4)]{kn}.

\begin{theorem}
\label{thm:rose} Let $S$ be an admissible set of paths in the rose graph
$R_{n}$, $n\geq2$, and let $\varepsilon=\sum_{k=1}^{n}e_{k}$.

\begin{enumerate}
\item[(a)] Let $c\in\mathrm{SCP}(R_{n})$ and $P=(p_{ij})\in GL_{n}(K)$. Then
$L\mathbf{1}_{S}\cong V_{[c^{\infty}]}^{P}$ if and only if $c=e_{j}$ for some
$j$ and the $j$-th column of $P$ is $(1,1,\dots,1)^{T}$. In this case
$\varphi_{P}(e_{j})=\varepsilon$ and $(L\mathbf{1}_{S})^{\varphi_{P}}\cong
V_{[e_{j}^{\infty}]}$.

\item[(b)] In particular, $L\mathbf{1}_{S}$ is isomorphic to a module of type
\textup{(5)} in \cite[Theorem 4.7]{nsv}; for instance $L\mathbf{1}_{S}\cong
V_{[e_{1}^{\infty}]}^{P}$, where $P\neq I_{n}$ is obtained from the identity
matrix $I_{n}$ by replacing its first column by $(1,1,\dots,1)^{T}$.

\item[(c)] $L\mathbf{1}_{S}$ is not isomorphic to any module of type
\textup{(1)}, \textup{(2)}, \textup{(3)} or \textup{(4)} in \cite[Theorem
4.7]{nsv}, nor to any module $V_{[c^{\infty}]}^{P}$ of type \textup{(5)} with
$(c,P)$ not as in \textup{(a)}.
\end{enumerate}
\end{theorem}

\begin{proof}
Recall that $L\mathbf{1}_{S}=(N_{S},\circ)\cong L/L(v-\varepsilon)$ by
Proposition~\ref{prop:rose}(b), and that $S=vS$ is infinite.

(a) Suppose $c=e_{j}$ and $p_{kj}=1$ for all $k$. Then $\varphi_{P}%
(e_{j})=\sum_{k}e_{k}p_{kj}=\varepsilon$, and by \eqref{eq:nsvc} and
Proposition~\ref{prop:rose}(b),
\[
V_{[e_{j}^{\infty}]}^{P}\cong L/L(v-\varepsilon)\cong L\mathbf{1}_{S}.
\]
Twisting by $\varphi_{P}$ gives $(L\mathbf{1}_{S})^{\varphi_{P}}%
\cong(V_{[e_{j}^{\infty}]}^{P})^{\varphi_{P}}=V_{[e_{j}^{\infty}]}$.

Conversely, suppose that $L\mathbf{1}_{S}\cong V_{[c^{\infty}]}^{P}$ and put
$q=\varphi_{P}(c)$ and $m=|c|\geq1$. Writing $c=e_{i_{1}}\cdots e_{i_{m}}$ and
expanding $q=\varphi_{P}(e_{i_{1}})\cdots\varphi_{P}(e_{i_{m}})$, we get
\begin{equation}
q=\sum_{|\mu|=m}a_{\mu}\mu,\qquad a_{\mu}=p_{k_{1}i_{1}}p_{k_{2}i_{2}}\cdots
p_{k_{m}i_{m}}\quad\text{for}\quad\mu=e_{k_{1}}e_{k_{2}}\cdots e_{k_{m}}.
\label{eq:qmu}%
\end{equation}
By \eqref{eq:nsvc}, the image $x\in N_{S}$ of the coset $v+L(v-q)$ under the
isomorphism is non-zero and satisfies $q\circ x=x$, that is, $x\cdot q^{\ast
}=x$. By \eqref{eq:right}, $(x\cdot\mu^{\ast})(p)=x(p^{\prime})$ if $p=\mu
p^{\prime}$, and $0$ otherwise. Since every $p\in D$ has a unique initial
segment $\mu_{p}$ of length $m$, we obtain
\begin{equation}
x(\mu_{p}p^{\prime})=a_{\mu_{p}}\,x(p^{\prime})\qquad(p=\mu_{p}p^{\prime}\in
D). \label{eq:xmu}%
\end{equation}
Iterating \eqref{eq:xmu}, for every finite path $\gamma$ of length $km$ there
is a scalar $A(\gamma)$, a product of $k$ of the coefficients $a_{\mu}$, such
that $x(\gamma y)=A(\gamma)x(y)$ for all $y\in D$.

By Lemma~\ref{lem:key} there are a finite path $\gamma_{0}$ and $0\neq
c_{0}\in K$ with $x(\gamma_{0}y)=c_{0}$ for all $y\in S$. Extend $\gamma_{0}$
to a path $\gamma=\gamma_{0}\delta$ whose length is a multiple of $m $. For
$y\in S$ we have $\delta y\in S$, so $A(\gamma)x(y)=x(\gamma y)=c_{0}\neq0$.
Hence $A(\gamma)\neq0$ and $x(y)=c_{0}A(\gamma)^{-1}$ for all $y\in S$; since
$\operatorname{supp}(x)\subseteq S$, we get $x=\lambda\mathbf{1}_{S}$ with
$\lambda\neq0$. Now let $\mu$ be any path of length $m$ and $y\in S$. Then
$\mu y\in S$, and \eqref{eq:xmu} gives $\lambda=x(\mu y)=a_{\mu}\lambda$, so
$a_{\mu}=1$.

By \eqref{eq:qmu}, $p_{k_{1}i_{1}}\cdots p_{k_{m}i_{m}}=1$ for all
$k_{1},\dots,k_{m}\in\{1,\dots,n\}$. Fixing all $k_{s}$ with $s\neq t$ and
varying $k_{t}$, we see that the $i_{t}$-th column of $P$ is a non-zero
multiple of $(1,1,\dots,1)^{T}$, for $t=1,\dots,m$. If two different loops
occurred in $c$, the corresponding two columns of $P$ would be proportional,
contradicting the invertibility of $P$. So $c=e_{j}^{m}$ for some $j$, and, as
$c$ is a simple closed path, $m=1$ and $c=e_{j}$. Finally $p_{kj}=a_{e_{k}}=1$
for all $k$.

(b) This follows from (a); note that $P\neq I_{n}$ since $n\geq2$.

(c) The module of type (1) is a Chen module of type (C1), so it is not
isomorphic to $L\mathbf{1}_{S}$ by Theorem~\ref{thm:chen}. If $L\mathbf{1}%
_{S}\cong S_{\Pi_{c}}^{f}\cong L/Lf(c)$, then the image $x_{0}\in N_{S}$ of
the coset of $v$ is non-zero and satisfies $v\circ x_{0}=x_{0}$ and $f(c)\circ
x_{0}=0$, where $f(0)=1$; this contradicts Proposition~\ref{prop:poly}, which
holds for every closed path $c$. By (b) and \cite[Theorem 4.5]{nsv},
$\operatorname{End}_{L}(L\mathbf{1}_{S})\cong K$, while $\operatorname{End}%
_{L}(S_{d}^{f,p})\cong K[x]/K[x]f(x)$ has dimension $\deg f\geq2$ over $K$; so
$L\mathbf{1}_{S}$ is not isomorphic to a module of type (3). By
Proposition~\ref{prop:rose}(b), $L\mathbf{1}_{S}$ is finitely presented,
whereas the modules of type (4) are not. The statement about type (5) is part (a).
\end{proof}

\begin{remark}
(1) By Theorem~\ref{thm:chen}, the simple module $L\mathbf{1}_{S}$ over
$L_{K}(R_{n})$ is not a Chen module, but by Theorem~\ref{thm:rose}(a) it
becomes the Chen module $V_{[e_{j}^{\infty}]}$ after twisting by the graded
automorphism $\varphi_{P}$, with $\mathbf{1}_{S}$ corresponding to
$e_{j}^{\infty}$: indeed $e_{j}\cdot\mathbf{1}_{S}=\varepsilon\circ
\mathbf{1}_{S}=\sum_{k}\mathbf{1}_{e_{k}S}=\mathbf{1}_{S}$ in $(L\mathbf{1}%
_{S})^{\varphi_{P}}$.

(2) For the right module itself, Proposition~\ref{prop:rose}(a) gives
$N_{S}\cong L/(v-\sum_{k}e_{k}^{\ast})L$, which is the right module obtained
from $V_{[e_{j}^{\infty}]}^{P}\cong L/L(v-\varepsilon)$ through the involution
of $L$.

(3) Since all the modules $N_{S}$ over $L_{K}(R_{n})$ are isomorphic,
Theorem~\ref{thm:rose} applies equally to a single tail-equivalence class
$S=[p]$, to a union of several classes, and to $S=D_{\infty}$, in which case
$N_{S}$ is the space of locally constant functions on $D_{\infty}$
(Example~\ref{ex:rose}).
\end{remark}

\section{The case of Cohn path algebras}

\label{sec:cohn}

In this section we explore the module generated by an indicator function in
the context of the Cohn path algebra $C:=C_{K}(E)$ instead of the Leavitt path
algebra $L_{K}(E)$. Recall \cite[Section 1.5]{aas} that $C$ is the $K$-algebra
generated by a set $\{v:v\in E^{0}\}$ of pairwise orthogonal idempotents
together with a set $\{e,e^{\ast}:e\in E^{1}\}$ satisfying the relations (1),
(2) and (3) (CK-1) of Definition~\ref{def:LPA}, but \emph{not} the relation
(CK-2). Thus $L_{K}(E)$ is the quotient of $C$ by the ideal generated by the
elements
\[
q_{u}=u-\sum_{e\in s^{-1}(u)}ee^{\ast}\qquad(u\text{ a regular vertex}),
\]
and we write $\pi:C\rightarrow L_{K}(E)$ for the quotient map. As in the case
of $L_{K}(E)$, the algebra $C$ is spanned by the monomials $\alpha\beta^{\ast
}$ with $\alpha,\beta$ finite paths and $r(\alpha)=r(\beta)$, and it is
$\mathbb{Z}$-graded, the monomial $\alpha\beta^{\ast}$ having degree
$|\alpha|-|\beta|$.

Since (CK-2) does not hold in $C$, every vertex behaves like a singular
vertex, and finite paths ending at regular vertices must also be allowed. So
the natural analogue of $D$ is the set
\[
D^{C}=D_{\infty}\cup\mathrm{Path}(E),
\]
where $\mathrm{Path}(E)$ denotes the set of \emph{all} finite paths in $E$,
including the vertices. Let $KD^{C}$ be the $K$-vector space with basis
$D^{C}$, and define $v\cdot p$, $e\cdot p$ and $e^{\ast}\cdot p$ for $p\in
D^{C}$ by the same rules as for $KD$ in Section~2.

\begin{lemma}
\label{lem:CKD}

\begin{enumerate}
\item[(a)] These rules make $KD^{C}$ a left $C$-module, and formula
\eqref{eq:action} holds for all $p\in D^{C}$.

\item[(b)] $KD$ is a $C$-submodule of $KD^{C}$, and the action of $C$ on $KD$
is the action of $L_{K}(E)$ given in Lemma~\ref{lem:KD}, composed with $\pi$.
\end{enumerate}
\end{lemma}

\begin{proof}
(a) The proof of Lemma~\ref{lem:KD} shows that the operators $T_{v}$, $T_{e}$,
$T_{e^{\ast}}$ on $KD^{C}$ satisfy the relations (1), (2) and (CK-1). The
assumption that the paths end at singular vertices was used there only to
verify (CK-2), which is not a relation of $C$. By the universal property of
$C$, $KD^{C}$ is a left $C$-module, and \eqref{eq:action} follows by applying
the generators one at a time.

(b) If $p\in D$, then $v\cdot p$, $e\cdot p$ and $e^{\ast}\cdot p$ belong to
$KD$, since a path obtained from a boundary path by adding or deleting an
initial segment is again a boundary path. The restriction of the action of $C$
to $KD$ and the composite of $\pi$ with the representation of
Lemma~\ref{lem:KD} are two homomorphisms $C\rightarrow\operatorname{End}%
_{K}(KD)$ which agree on the generators of $C$, so they coincide.
\end{proof}

As in Section~2, $K^{D^{C}}=\operatorname{Hom}_{K}(KD^{C},K)$ is a right
$C$-module under $(f\cdot a)(p)=f(a\cdot p)$, and formula \eqref{eq:right}
holds for all $f\in K^{D^{C}}$ and $p\in D^{C}$. We use the notions $(\#)$,
$\beta\preceq p$, $S^{0}$, $vS$, $\beta S$, $s_{S}(v)$, strongly connected and
admissible exactly as in Definition~\ref{def:adm}, with $D$ replaced by
$D^{C}$. For $S\subseteq D^{C}$ we put
\[
N_{S}^{C}=\mathbf{1}_{S}\cdot C\subseteq K^{D^{C}}\qquad\text{and}\qquad
H=E^{0}\setminus S^{0}.
\]

\begin{lemma}
\label{lem:Cind} Let $S\subseteq D^{C}$ satisfy $(\#)$. Then $\mathbf{1}%
_{S}\cdot\alpha\beta^{\ast}=\mathbf{1}_{\beta S}$ for all finite paths
$\alpha,\beta$ with $r(\alpha)=r(\beta)$. Consequently, the statements of
Lemmas~\ref{lem:ind} and \ref{lem:local} hold for $N_{S}^{C}$.
\end{lemma}

\begin{proof}
The proofs of Lemmas~\ref{lem:ind} and \ref{lem:local} use only formula
\eqref{eq:right}, the property $(\#)$, and the fact that the monomials
$\alpha\beta^{\ast}$ span the algebra. All three are available for $C$.
\end{proof}

\begin{lemma}
\label{lem:transfer} Let $S\subseteq D$ satisfy $(\#)$. Then the restriction
map $\rho:K^{D^{C}}\rightarrow K^{D}$, $f\mapsto f|_{D}$, is a homomorphism of
right $C$-modules which maps $N_{S}^{C}$ bijectively onto $N_{S}%
=\mathbf{1}_{S}\cdot L_{K}(E)$, where $C$ acts on $N_{S}$ through $\pi$. In
particular, $N_{S}^{C}$ and $N_{S}$ have the same submodules.
\end{lemma}

\begin{proof}
Let $f\in K^{D^{C}}$, $a\in C$ and $p\in D$. Then $a\cdot p\in KD$ by
Lemma~\ref{lem:CKD}(b), so $(f\cdot a)(p)=\rho(f)(a\cdot p)=(\rho(f)\cdot
a)(p)$, and $\rho$ is a homomorphism. By Lemma~\ref{lem:Cind}, every element
of $N_{S}^{C}$ is supported in $S\subseteq D$, so $\rho$ is injective on
$N_{S}^{C}$. Finally, $\rho(\mathbf{1}_{S}\cdot a)=\mathbf{1}_{S}\cdot\pi(a)$
for $a\in C$, and $\pi$ is surjective, so $\rho(N_{S}^{C})=N_{S}$.
\end{proof}

\begin{lemma}
\label{lem:Cadm} Let $S\subseteq D^{C}$.

\begin{enumerate}
\item[(a)] The statements and proofs of Lemma~\ref{lem:basic}(a), (b), (c),
(e) and of Lemma~\ref{lem:single} hold for $S$. In particular, $H=E^{0}%
\setminus S^{0}$ is hereditary.

\item[(b)] The set $\{w\}$ is admissible if and only if $w$ is a source; here
$w$ may be a sink, a regular vertex or an infinite emitter.

\item[(c)] If $S\subseteq D$ satisfies $(\#)$, then $H$ is also saturated.
\end{enumerate}
\end{lemma}

\begin{proof}
(a) These proofs use only $(\#)$, strong connectivity and the combinatorics of
paths; they never use that the finite paths in $S$ end at singular vertices.

(b) If $w$ is a source, then the only finite path ending at $w$ is $w$ itself,
so $\{w\}$ satisfies $(\#)$, and it is clearly strongly connected. Conversely,
if $\{w\}$ is admissible and $e$ is an edge with $r(e)=w$, then $e=e\,w\in
\{w\}$ by $(\#)$, which is absurd.

(c) The proof of Lemma~\ref{lem:basic}(d) applies, since a regular vertex does
not belong to $D$.
\end{proof}

In general $H$ need not be saturated. For example, if $E$ consists of a single
edge from a vertex $w$ to a sink $u$, then $S=\{w\}$ is admissible in $D^{C}$
by Lemma~\ref{lem:Cadm}(b), and $H=\{u\}$ is not saturated. So the breaking
vertices of $H$ are not available in general, and we replace the condition of
Theorem~\ref{thm:simple} by the following one.

\begin{definition}
An admissible set $S\subseteq D^{C}$ satisfies \emph{condition (B)} if either
$S=\{w\}$ for a vertex $w$, or $s_{S}(w)$ is infinite for every vertex $w$
which is an element of $S$.
\end{definition}

\begin{lemma}
\label{lem:Bequiv} Let $S\subseteq D$ be admissible. Then $S$ satisfies
\textup{(B)} if and only if $S$ contains no breaking vertex of $H$.
\end{lemma}

\begin{proof}
By Lemma~\ref{lem:Cadm}(c), $H$ is hereditary and saturated, so $B_{H}$ is
defined. If $S=\{w\}$, then $s_{S}(w)=\emptyset$, so $w\notin B_{H}$. Suppose
that $S$ is not a single vertex, and let $w$ be a vertex which is an element
of $S$. Then $w$ is singular, since $S\subseteq D$, and $s_{S}(w)\neq
\emptyset$ by Lemma~\ref{lem:single}; so $w$ is an infinite emitter. As
$s_{S}(w)=s^{-1}(w)\cap r^{-1}(E^{0}\setminus H)$, we have $w\in B_{H}$ if and
only if $s_{S}(w)$ is finite.
\end{proof}

\begin{theorem}
\label{thm:Csimple} Let $S\subseteq D^{C}$ be an admissible set. Then
$N_{S}^{C}=\mathbf{1}_{S}\cdot C$ is a simple right $C$-module if and only if
$S$ satisfies \textup{(B)}. Moreover, if $S$ satisfies \textup{(B)}, then
either $S\subseteq D$, or $S=\{w\}$ for a regular vertex $w$ which is a source.
\end{theorem}

\begin{proof}
Suppose first that $S=\{w\}$. Then $w$ is a source by Lemma~\ref{lem:Cadm}(b),
so the only finite path ending at $w$ is $w$ itself, and $\mathbf{1}%
_{S}=\delta_{w}$. By Lemma~\ref{lem:Cind}, $N_{S}^{C}$ is spanned by the
functions $\mathbf{1}_{\beta S}$ with $r(\beta)\in S^{0}=\{w\}$; the only such
$\beta$ is $w$, and $\mathbf{1}_{wS}=\delta_{w}$. Hence $N_{S}^{C}=K\delta
_{w}$ is one-dimensional over $K$. If $U$ is a non-zero submodule, then
$c\delta_{w}\in U$ for some $0\neq c\in K$, and $\delta_{w}=(c\delta_{w}%
)\cdot(c^{-1}w)\in U$; so $U=N_{S}^{C}$, and $N_{S}^{C}$ is simple.

Assume now that $S$ is not a single vertex. Suppose that (B) holds. Then the
proof of Lemma~\ref{lem:key} applies verbatim: its Case~2 uses only that
$s_{S}(w)$ is infinite for the vertex $w=r(p)$, which belongs to $S$, and this
is exactly (B). The first half of the proof of Theorem~\ref{thm:simple} then
shows that $g\cdot C=N_{S}^{C}$ for every non-zero $g\in N_{S}^{C}$, so
$N_{S}^{C}$ is simple.

Conversely, suppose that $w$ is a vertex which is an element of $S$ such that
$s_{S}(w)$ is finite. By Lemma~\ref{lem:single}, $s_{S}(w)=\{e_{1},\dots
,e_{k}\}$ with $k\geq1$. The second half of the proof of
Theorem~\ref{thm:simple} applies word for word (it does not use that $w$ is an
infinite emitter), and gives
\[
\delta_{w}=\mathbf{1}_{S}\cdot\Bigl(w-\sum_{i=1}^{k}e_{i}e_{i}^{\ast}\Bigr)\in
N_{S}^{C}\qquad\text{and}\qquad\delta_{w}\cdot C=\operatorname{span}%
_{K}\{\delta_{\beta}:r(\beta)=w\},
\]
where $\beta$ runs over the finite paths ending at $w$. The latter is a
non-zero submodule consisting of functions with finite support, and it does
not contain $\mathbf{1}_{wS}$, which has infinite support. So $N_{S}^{C}$ is
not simple.

For the last statement, suppose that (B) holds and that $S$ is not a single
vertex. If $p\in S$ is a finite path, then $r(p)$ is a vertex belonging to $S$
(Lemma~\ref{lem:basic}(a)), so $s_{S}(r(p))$ is infinite, $r(p)$ is an
infinite emitter, and $p\in D$. Hence $S\subseteq D$. If $S=\{w\}$ and $w$ is
not regular, then $w$ is singular and again $S\subseteq D$.
\end{proof}

\begin{corollary}
\label{cor:classify} Up to isomorphism, the simple modules $N_{S}^{C}$ are the following.

\begin{enumerate}
\item[(i)] The simple $L_{K}(E)$-modules $N_{S}$, where $S\subseteq D$ is
admissible and contains no breaking vertex of $H$, regarded as $C$-modules
through $\pi$; here $N_{S}^{C}\cong N_{S}$ by Lemma~\ref{lem:transfer}.

\item[(ii)] For each regular vertex $w$ which is a source, the module
$N_{\{w\}}^{C}=K\delta_{w}$. It is isomorphic to $q_{w}C=Kq_{w}$, and so it is
one-dimensional, graded, projective and finitely presented. It is not
annihilated by $\ker\pi$; hence it is not an $L_{K}(E)$-module, and it is not
isomorphic to any of the modules in \textup{(i)}.
\end{enumerate}
\end{corollary}

\begin{proof}
By Theorem~\ref{thm:Csimple}, a simple module $N_{S}^{C}$ either has
$S\subseteq D$, and then (i) follows from Lemmas~\ref{lem:transfer} and
\ref{lem:Bequiv} and Theorem~\ref{thm:simple}, or $S=\{w\}$ with $w$ a regular source.

Let $w$ be a regular source. The idempotents $ee^{\ast}$, $e\in s^{-1}(w)$,
are pairwise orthogonal by (CK-1) and satisfy $w\,ee^{\ast}=ee^{\ast}%
=ee^{\ast}w$, so $q_{w}$ is an idempotent of degree $0$. We claim that
$q_{w}C=Kq_{w}$. If $f\in s^{-1}(w)$, then $q_{w}f=f-\sum_{e\in s^{-1}%
(w)}ee^{\ast}f=f-f=0$ by (CK-1). Hence, for a monomial $\alpha\beta^{\ast}$,
we have $q_{w}\alpha\beta^{\ast}=0$ unless $\alpha=w$; and if $\alpha=w$, then
$r(\beta)=w$ forces $\beta=w$, as $w$ is a source. This proves the claim.

By \eqref{eq:right}, $\delta_{w}\cdot ee^{\ast}=0$ for every edge $e$, since
$ep^{\prime}\neq w$ for every path $p^{\prime}$. Hence $\delta_{w}\cdot
q_{w}=\delta_{w}\cdot w=\delta_{w}$, and $a\mapsto\delta_{w}\cdot a$ is a
non-zero, hence bijective, homomorphism $Kq_{w}\rightarrow K\delta_{w}$. Thus
$K\delta_{w}\cong q_{w}C$, which is a graded direct summand of $C_{C}$; in
particular, it is projective and finitely presented. Since $q_{w}\in\ker\pi$
and $\delta_{w}\cdot q_{w}\neq0$, the module $K\delta_{w}$ is not an
$L_{K}(E)$-module, whereas the modules in (i) are annihilated by $\ker\pi$.
\end{proof}

\begin{corollary}
\label{cor:Cinf}

\begin{enumerate}
\item[(a)] If $S\subseteq D_{\infty}$ is an admissible set of infinite paths,
then $N_{S}^{C}$ is a simple $C$-module.

\item[(b)] If $E$ is a row-finite graph and $S\subseteq D^{C}$ is admissible,
then $N_{S}^{C}$ is simple if and only if either $S\subseteq D_{\infty}$, or
$S=\{w\}$ for a vertex $w$ which is a source.
\end{enumerate}
\end{corollary}

\begin{proof}
(a) $S$ contains no vertex, so condition (B) holds vacuously.

(b) In a row-finite graph, $s_{S}(w)$ is finite for every vertex $w$. So (B)
says that either $S$ is a single vertex, or $S$ contains no vertex, that is,
by Remark~\ref{important}(3), no finite path.
\end{proof}

Thus Corollary~\ref{cor:infinite} remains valid for Cohn path algebras, while
Corollary~\ref{Finite case} does not. For instance, if $E$ has one vertex $v$
and one loop $c$, then $S=\{c^{n}:n\geq0\}$ is an admissible subset of $D^{C}$
which is not a single vertex, so $N_{S}^{C}$ is not simple by
Corollary~\ref{cor:Cinf}(b).

\subsection*{Acknowledgement}

The authors thank Zak Mesyan for many helpful conversations and suggestions
that improved this paper considerably. The authors also used the AI assistant
Claude Opus 5.5 (Anthropic) to help improve the quality and presentation of
the results in this paper. We take complete responsibility for the contents of
this paper.


\begin{thebibliography}{99}                                                                                               %


\bibitem {aas}G. Abrams, P. Ara and M. Siles Molina, \emph{Leavitt path
algebras}, Lecture Notes in Mathematics, vol.~2191, Springer, 2017.

\bibitem {an}P.N. \'Anh and T.G. Nam, Special irreducible representations of
Leavitt path algebras, Adv. Math. 377 (2021), 107483.

\bibitem {ar}P. Ara and K.M. Rangaswamy, Finitely presented simple modules
over Leavitt path algebras, J. Algebra 417 (2014), 333--352.

\bibitem {chen}X.W. Chen, Irreducible representations of Leavitt path
algebras, Forum Math. 27 (2015), 549--574.

\bibitem {hr}R. Hazrat and K.M. Rangaswamy, On graded irreducible
representations of Leavitt path algebras, J. Algebra 450 (2016), 458--486.

\bibitem {kn}S. Kuroda and T.G. Nam, Anick type automorphisms and new
irreducible representations of Leavitt path algebras, J. Noncommut. Geom. 17
(2023), 811--834.

\bibitem {nsv}T.G. Nam, A.K. Srivastava and N.T. Vien, Automorphisms of
Leavitt path algebras: Zhang twist and irreducible representations, J. Algebra
654 (2024), 189--234.

\bibitem {pacheco}E. Pacheco, \emph{New simple representations of Leavitt path
algebras}, Ph.D. thesis, University of Oklahoma, 2019.

\bibitem {r-prime}K.M. Rangaswamy, The theory of prime ideals of Leavitt path
algebras over arbitrary graphs, J. Algebra 375 (2013), 73--96.

\bibitem {r-simple}K.M. Rangaswamy, On simple modules over Leavitt path
algebras, J. Algebra 423 (2015), 239--258.

\bibitem {T}M. Tomforde, Uniqueness theorems and ideal structure for Leavitt
path algebras, J. Algebra 318 (2007), 270--299.
\end{thebibliography}
\end{document}